\documentclass[11pt]{article} \usepackage{amsfonts} \usepackage{latexsym} \usepackage{graphicx}  

\newcommand{\breath}{\vspace{3mm}} 
\newcommand{\pause}{\breath                     \hrulefill }    
\newcounter{thmcount} 

\newcounter{excount}  

\newcounter{claimcount}[thmcount] 

\newcounter{subclaimcount}[claimcount] 

\newcounter{subsubclaimcount}[subclaimcount]  

\newsavebox{\thmstyle}  

\newfont{\fibfont}{cmfib8}    

\newcommand{\thmfont}[1]{{\sl #1}} 
\newcommand{\Theorem}[2]{        \refstepcounter{thmcount}        \subparagraph{{\normalsize Theorem  \thethmcount:}} {\em #1} 			 \savebox{\thmstyle}{{\tt Theorem}}                            \thmfont{ #2 }}
\newcommand{\Corollary}[2]{        \refstepcounter{thmcount}        \subparagraph{{\normalsize Corollary  \thethmcount:}} {\em #1}  	  		 \savebox{\thmstyle}{{\tt Corollary}}                            \thmfont{ #2 }}
\newcommand{\Proposition}[2]{        \refstepcounter{thmcount}        \subparagraph{{\normalsize Proposition  \thethmcount:}} {\em #1} 			 \savebox{\thmstyle}{{\tt Proposition}}                            \thmfont{ #2 }}
\newcommand{\Definition}[2]{        \refstepcounter{thmcount}        \subparagraph{{\normalsize Definition  \thethmcount:}} {\em #1} 			 \savebox{\thmstyle}{{\tt Definition }}                        {\begin{list}{} 			{\setlength{\leftmargin}{2em} 			\setlength{\rightmargin}{0em}}                         \item \thmfont{ #2 }  \end{list} }          \breath}

\newcommand{\Claim}[1]{\refstepcounter{claimcount}                {\bf Claim \theclaimcount: \ }\thmfont{ #1}} 
\newcommand{\claim}{\Claim}
\newcommand{\thmpart}[1]{{\bf Part #1}} 
\newcommand{\Proofof}[1]{{\bf \hspace{-1em}  Proof of #1: \ \ }}
\newcommand{\proofof}[1]{\Proofof{#1}} 
\newenvironment{prf}[1]{\begin{list}{} 			{\setlength{\leftmargin}{1em} 			\setlength{\rightmargin}{0em}}                         \item {\bf \hspace{-1em}  #1 \ \ }}                     {\end{list}  }   

\newenvironment{thmproof}{ \begin{prf}{Proof:}}                      {                {\tt \hrulefill $\Box$} \end{prf}                  \breath  }  

\newenvironment{claimproof}{ \begin{prf}{Proof:}}                      {                {\tt \dotfill $\Box$[Claim \theclaimcount]} \end{prf}}

\newcommand{\Proof}{     {\bf Proof:}}
 
\newcommand{\QED}{\hrulefill\ensuremath{\Box}}
\newcommand{\qed}{\QED}

\newcommand{\btxt}[1]{\mbox{#1}}
\newcommand{\etxt}{}                           
\newcommand{\beq}{\begin{eqnarray*}}
\newcommand{\eeq}{\end{eqnarray*}} 
\newcommand{\beqn}{\begin{equation}}
\newcommand{\eeqn}{\end{equation}} 
\newcommand{\blist}{\begin{enumerate}}
\newcommand{\elist}{\end{enumerate}} 
\newcommand{\bitem}{\begin{itemize}}
\newcommand{\eitem}{\end{itemize}} 
\newcommand{\bquote}{\begin{quotation}}
\newcommand{\equote}{\end{quotation}}    
\newcommand{\Cesaro}{Ces\`aro }
\newcommand{\done}{\ensuremath{\mathsf{ 1\!\!1}}} 
\newcommand{\dC}{\ensuremath{\mathbb{C}}}
\newcommand{\dG}{\ensuremath{\mathbb{G}}}
\newcommand{\dH}{\ensuremath{\mathbb{H}}}
\newcommand{\dM}{\ensuremath{\mathbb{M}}}
\newcommand{\dN}{\ensuremath{\mathbb{N}}}
\newcommand{\dP}{\ensuremath{\mathbb{P}}}
\newcommand{\dR}{\ensuremath{\mathbb{R}}}

\newcommand{\dT}{\ensuremath{\mathbb{T}}}

\newcommand{\dZ}{\ensuremath{\mathbb{Z}}}

\newcommand{\barchi}{\ensuremath{{\bar{\chi }}}}
\newcommand{\bC}{\ensuremath{\mathbf{ C}}}
\newcommand{\bL}{\ensuremath{\mathbf{ L}}}
\newcommand{\bN}{\ensuremath{\mathbf{ N}}}
\newcommand{\bP}{\ensuremath{\mathbf{ P}}}
\newcommand{\bQ}{\ensuremath{\mathbf{ Q}}}
\newcommand{\bR}{\ensuremath{\mathbf{ R}}}
\newcommand{\bS}{\ensuremath{\mathbf{ S}}}
\newcommand{\bT}{\ensuremath{\mathbf{ T}}}
\newcommand{\bU}{\ensuremath{\mathbf{ U}}}
\newcommand{\bZ}{\ensuremath{\mathbf{ Z}}} 
\newcommand{\ba}{\ensuremath{\mathbf{ a}}}
\newcommand{\bb}{\ensuremath{\mathbf{ b}}}
\newcommand{\bc}{\ensuremath{\mathbf{ c}}}

\newcommand{\bg}{\ensuremath{\mathbf{ g}}}
\newcommand{\bh}{\ensuremath{\mathbf{ h}}}
\newcommand{\bi}{\ensuremath{\mathbf{ i}}}
\newcommand{\bk}{\ensuremath{\mathbf{ k}}}

\newcommand{\bn}{\ensuremath{\mathbf{ n}}}

\newcommand{\bp}{\ensuremath{\mathbf{ p}}}

\newcommand{\br}{\ensuremath{\mathbf{ r}}}
\newcommand{\bt}{\ensuremath{\mathbf{ t}}}
\newcommand{\bu}{\ensuremath{\mathbf{ u}}}
\newcommand{\bv}{\ensuremath{\mathbf{ v}}}
\newcommand{\bw}{\ensuremath{\mathbf{ w}}}

\newcommand{\sA}{\ensuremath{\mathcal{ A}}}
\newcommand{\sB}{\ensuremath{\mathcal{ B}}}
\newcommand{\sC}{\ensuremath{\mathcal{ C}}}
\newcommand{\sE}{\ensuremath{\mathcal{ E}}}
\newcommand{\sF}{\ensuremath{\mathcal{ F}}}

\newcommand{\sL}{\ensuremath{\mathcal{ L}}}
\newcommand{\sM}{\ensuremath{\mathcal{ M}}}
\newcommand{\sN}{\ensuremath{\mathcal{ N}}}
\newcommand{\sP}{\ensuremath{\mathcal{ P}}}
\newcommand{\sR}{\ensuremath{\mathcal{ R}}}
\newcommand{\sS}{\ensuremath{\mathcal{ S}}}

\newcommand{\sU}{\ensuremath{\mathcal{ U}}}
\newcommand{\sV}{\ensuremath{\mathcal{ V}}}
\newcommand{\sW}{\ensuremath{\mathcal{ W}}}
\newcommand{\sX}{\ensuremath{\mathcal{ X}}}

\newcommand{\gA}{\ensuremath{\mathfrak{ A}}}
\newcommand{\gB}{\ensuremath{\mathfrak{ B}}}

\newcommand{\gR}{\ensuremath{\mathfrak{ R}}}
\newcommand{\gT}{\ensuremath{\mathfrak{ T}}}

\newcommand{\gW}{\ensuremath{\mathfrak{ W}}}
\newcommand{\ga}{\ensuremath{\mathfrak{ a}}}

\newcommand{\alp }{\ensuremath{\alpha}}
\newcommand{\bet }{\ensuremath{\beta}}
\newcommand{\gam }{\ensuremath{\gamma}}
\newcommand{\del }{\ensuremath{\delta}}
\newcommand{\eps }{\ensuremath{\epsilon}}

\newcommand{\hT}{\ensuremath{{\widehat{T}}}}
\newcommand{\hX}{\ensuremath{{\widehat{X}}}}

\newcommand{\hsX}{\ensuremath{{\widehat{\mathcal{ X}}}}}
\newcommand{\hdel }{\ensuremath{\widehat{\delta}}}
\newcommand{\hmu }{\ensuremath{\widehat{\mu}}}
\newcommand{\hnu }{\ensuremath{\widehat{\nu}}}

\newcommand{\fA}{\ensuremath{\mathsf{ A}}}
\newcommand{\fB}{\ensuremath{\mathsf{ B}}}
\newcommand{\fC}{\ensuremath{\mathsf{ C}}}
\newcommand{\fa}{\ensuremath{\mathsf{ a}}}
\newcommand{\fb}{\ensuremath{\mathsf{ b}}}
\newcommand{\fc}{\ensuremath{\mathsf{ c}}}
\newcommand{\fd}{\ensuremath{\mathsf{ d}}}

\newcommand{\tlS}{\ensuremath{{\widetilde{S}}}}
\newcommand{\tla}{\ensuremath{{\widetilde{a}}}}

\newcommand{\tlc}{\ensuremath{{\widetilde{c}}}}

\newcommand{\tlm}{\ensuremath{{\widetilde{m}}}}

\newcommand{\tlu}{\ensuremath{{\widetilde{u}}}}
\newcommand{\tlbC}{\ensuremath{{\widetilde{\mathbf{ C}}}}}
\newcommand{\tlba}{\ensuremath{{\widetilde{\mathbf{ a}}}}}
\newcommand{\tlbc}{\ensuremath{{\widetilde{\mathbf{ c}}}}}
\newcommand{\tlbg}{\ensuremath{{\widetilde{\mathbf{ g}}}}}
\newcommand{\tlsM}{\ensuremath{{\widetilde{\mathcal{ M}}}}}
\newcommand{\tlsU}{\ensuremath{{\widetilde{\mathcal{ U}}}}}
\newcommand{\tlsV}{\ensuremath{{\widetilde{\mathcal{ V}}}}}
\newcommand{\tlmu }{\ensuremath{{\widetilde{\mu}}}}
\newcommand{\tlnu }{\ensuremath{{\widetilde{\nu}}}}

\newcommand{\lb}{\ensuremath{\left}}
\newcommand{\rb}{\ensuremath{\right}} 
 
\newcommand{\maketall}{\rule[-0.5cm]{0cm}{1cm}}       
\newcommand{\implies}{\mbox{$\Longrightarrow$}}
\newcommand{\map}{\ensuremath{\longrightarrow}}

\newcommand{\goto}{\ensuremath{\rightarrow}}
\newcommand{\into}{\ensuremath{\map}}
\newcommand{\seilpmi}{\ensuremath{\Longleftarrow}}
\newcommand{\statement}[1]{\lb( \ \maketall       \begin{minipage}{40em}       \begin{tabbing}         #1        \end{tabbing}      \end{minipage} \ \rb)}     
\newcommand{\oo}{\ensuremath{\infty}}        
\newcommand{\X}{\ensuremath{\times}}
\newcommand{\x}{\ensuremath{\X}}
\newcommand{\tensor}{\ensuremath{\otimes}}
\newcommand{\Tensor}{\ensuremath{\bigotimes}}
\newcommand{\dirsum}{\ensuremath{\oplus}}
\newcommand{\intsct}{\ensuremath{\cap}}
\newcommand{\Intsct}{\ensuremath{\bigcap}}
\newcommand{\disj}{\ensuremath{\sqcup}}
\newcommand{\Disj}{\ensuremath{\bigsqcup}}   
\newcommand{\set}[2]{\ensuremath{\left\{ #1 \; ; \; #2 \right\} }} 
\newcommand{\seq}[2]{\ensuremath{ \lb\{#1 |_{_{{#2}}} \rb\} }}     
\newcommand{\supp}[1]{ \ensuremath{\mathbf{ supp}\lb[#1\rb]}}    
\newcommand{\norm}[2]{\ensuremath{\left\| #1 \right\|_{{#2}} }   }
\newcommand{\inn}[1]{\ensuremath{\left\langle #1 \right\rangle }}       
\newcommand{\Id}[1]{\ensuremath{\mathbf{ Id}_{{#1}}}}
\newcommand{\pr}[1]{\ensuremath{\mathbf{ pr}_{{#1}}}}
 
\newcommand{\chr}[1]{\ensuremath{{\done}_{{#1}}}}

\newcommand{\interior}[1]{\ensuremath{\bi\bn\bt\lb[ #1 \rb]}}  
\newcommand{\shift}[1]{ {{\sS}_{^{\!\!h\!i\!f\!t}}^{#1}} }     
\newsavebox{\gaussopen} \savebox{\gaussopen}{\begin{picture}(5,5)                    \put(0,5){\oval(12,12)[br]}                    \put(12,4){\oval(12,12)[tl]}                      \end{picture}                   }  

\newsavebox{\gaussclose} \savebox{\gaussclose}{\begin{picture}(5,5)                    \put(9,5){\oval(12,12)[bl]}                    \put(-3,4){\oval(12,12)[tr]}                      \end{picture}                   }

\newcommand{\Lbsg}{{\ensuremath{\sL^{_{\!b\!s\!g}}}}}
\newcommand{\card}[1]{\ensuremath{{\sC_{^{\!\!a\!r\!d}}\lb[#1\rb]}}}
\newcommand{\Meas}[1]{\ensuremath{\sM_{^{\!\sE\!\!\sA\!\sS}} \lb[#1\rb] }}
\newcommand{\invMeas}[2]{\ensuremath{\sM^{#1}_{^{\!\sE\!\!\sA\!\sS}} \lb[#2\rb] }}
\newcommand{\marg}[2]{\ensuremath{\pr{#1}^* \lb[#2\rb]}}     
\newcommand{\Four}{\ensuremath{{\sF_{^{\!\!\!o\!u\!r}}}}}

\newcommand{\mtrx}[3]{\ensuremath{\lb[#1  |_{#2}^{#3} \rb]}}         
 
\newcommand{\Real}{\dR}
\newcommand{\Natur}{\dN}
\newcommand{\Zahl}{\dZ}
\newcommand{\Cplx}{\dC}
\newcommand{\Torus}[1]{\ensuremath{{\dT}^{#1}}}  
\newsavebox{\cubebox} \savebox{\cubebox}{\begin{picture}(12,8)                        \put(0,-5){\framebox(10,10)}                    \put(3,-2){\framebox(10,10)}                    \end{picture}                   } 

\newsavebox{\ballbox} \savebox{\ballbox}{\begin{picture}(9,8)                    \put(5,3){\circle{15}}                    \put(5,3){\circle{12}}                    \end{picture}                   }  

\newsavebox{\cycopen} \savebox{\cycopen}{\begin{picture}(5,5)                   \put(0,4){\oval(12,12)[l]}                                        \put(1,-2){\vector(1,0){5}}                    \end{picture}                   }  

\newsavebox{\cycclose} \savebox{\cycclose}{\begin{picture}(5,5)                   \put(9,4){\oval(12,12)[r]}                                        \put(7,10){\vector(-1,0){5}}                    \end{picture}                    }  

\newcommand{\CC}[1]{\ensuremath{{\lb[ #1 \rb]}}}
\newcommand{\CO}[1]{\ensuremath{{\lb[ #1 \rb)}}}
\newcommand{\pb}[1]{#1^{\!\nwarrow}}

\newcommand{\tlgW}{\widetilde{\gW}} 
\newcommand{\Latt}{{\Zahl^D}}
\newcommand{\statMeas}[1]{\invMeas{s\!t\!a\!t}{#1}}
\newcommand{\extMeas}[1]{\invMeas{e\!x\!t}{#1}}
\newcommand{\ergMeas}[1]{\invMeas{e\!r\!g}{#1}} 
\newcommand{\dual}[1]{\widehat{\gA^{#1}}} 
\newcommand{\hgA}{\widehat{\gA}} 
\newcommand{\Freq}[3]{{\bf Freq}\lb[ #1 \subset #2; \ #3\rb] } 
\newcommand{\Fin}[1]{\sF_{\!i\!n}\lb[#1\rb]}   
\newcommand{\tilebox}[2]{\fbox{$\begin{array}{#1} #2 \end{array}$}}  
\newcommand{\connect}{\!\!\!\!\!\longleftrightarrow\!\!\!\!\!\!}
\newcommand{\minitilebox}[4]{\fbox{${^{#1}_{#3} \; ^{#2}_{#4} }$}} 

\newcommand{\RReal}{\Real_\dagger}
\newcommand{\RMeas}[1]{{\sM_{^{\!\sE\!\!\sA\!\sS}\dagger} \lb[#1\rb] }}
\newcommand{\RinvMeas}[2]{{\sM^{#1}_{^{\!\sE\!\!\sA\!\sS}\dagger} \lb[#2\rb] }}
\newcommand{\RstatMeas}[1]{\RinvMeas{s\!t\!a\!t}{#1}}
\newcommand{\RextMeas}[1]{\RinvMeas{e\!x\!t}{#1}}

\begin{document}

\title{Building a Stationary Stochastic Process from a Finite-dimensional
Marginal}

\author{Marcus Pivato}

\maketitle

\abstract{If $\gA$ is a finite alphabet, $\sU \subset \Latt$, 
and $\mu_\sU$ is a probability measure on $\gA^\sU$ that ``looks like''
the marginal projection of a stationary stochastic process on $\gA^\Latt$,
then can we ``extend'' $\mu_\sU$ to such a process?  Under what conditions
can we make this extension ergodic, (quasi)periodic, or (weakly) mixing?
After surveying classical work on this problem when $D = 1$,
we provide some sufficient conditions and some necessary conditions for
$\mu_\sU$ to be extendible for $D > 1$, and show that, in general, the 
problem is not formally decidable.}  

\subsubsection*{Mathematics Subject Classification Number:}

\subparagraph{Primary:} 37A50, 60G10 (Ergodic Theory of Stationary Stochastic Processes)
\subparagraph{Secondary:} 37B10 (Symbolic Dynamics)

\newpage

\section{Introduction}

\subsection{The Markov Extension in $\Zahl$} 

  Let $\gA$ be a finite alphabet, and let $\gA^\Zahl$ be the space of
{\em bi-infinite sequences} on $\gA$.  A {\bf stationary stochastic
process} is a probability measure $\mu$ on $\gA^\Zahl$ so that,
for any $V \in \Natur, \ \ b_0,b_1,\ldots,b_V \in \gA$, and any $k \in \Zahl$

\beq
\lefteqn{
 \mu\set{\ba \in \gA^\Zahl}{a_0 = b_0,\ldots,a_V = b_V }
}\\
  &\hspace{5em} = & 
  \mu\set{\ba \in \gA^\Zahl}{a_k = b_0,\ldots,a_{k+V} = b_{V} } \\
\eeq

  Let $\sU$ be the interval $\CC{0...U} \subset \Zahl$.  The {\bf projection
map} $\pr{\sU} : \gA^\Zahl \into  \gA^\sU$ is the map sending the
sequence $\mtrx{a_n}{n \in \Zahl}{}$ to the sequence $\mtrx{a_n}{n \in \sU}{}$.
With this map, we can project $\mu$ down to a {\bf marginal} measure,
$\mu_\sU \ := \ \marg{\sU} \mu$, on the space $\gA^\sU$.  This marginal
is then {\bf locally stationary}:  for any $V < U$, any
 $b_0,b_1,\ldots,b_V \in \gA$, and any $k \in \Zahl$ so that $V + k \leq U$
also,

\beq
 \lefteqn{
\mu_\sU\set{\ba \in \gA^\sU}{a_0 = b_0,\ldots,a_V = b_V } 
}\\
  &\hspace{5em} = & 
  \mu_\sU\set{\ba \in \gA^\sU}{a_k = b_0,\ldots,a_{k+V} = b_{V} } \\
\eeq
               
  Can we {\em reverse} this process?  Given a locally stationary
measure $\mu_\sU$ upon $\gA^\sU$, can we {\bf extend} it to a stationary
stochastic process $\mu$ on $\gA^\Zahl$, so that
$\marg{\sU} \mu \ = \ \mu_\sU$?  Yes, and furthermore, we can do so
in a canonical fashion, via the so-called {\bf Markov Extension}. 

  An intuitive description of the Markov Extension is this:
We randomly ``choose'' the coordinates $a_0,\ldots,a_U$ according
to the probability measure $\mu_\sU$.  We then randomly chose the
coordinate $a_{U+1}$, again according to $\mu_\sU$ (now treated
as a probability measure on $\gA^{\sU+1}$), but {\em conditioned}
upon the fact that we have already fixed coordinates $a_1,\ldots,a_U$.
Next, we randomly chose the
coordinate $a_{U+2}$, again according to $\mu_\sU$ (now treated
as a probability measure on $\gA^{\sU+2}$), but {\em conditioned}
upon the fact that we have already fixed coordinates $a_2,\ldots,a_{U+1}$.
Inductively, we get a $U-$step Markov process
on $\gA$.

 To formally construct the Markov Extension, we need a bit
of notation:

\bitem
  \item  If $\ba \ = \ \mtrx{a_n}{n\in\Zahl}{}$ is an element of $\gA^\Zahl$,
and $\sV \subset \Zahl$, then let $\ba_\sV \ := \ \mtrx{a_v}{v \in \sV}{}$
   
  \item  If $\mu$ is a measure upon $\gA^\Zahl, \ \sV \subset \Zahl$, 
and $\bb \in \gA^\sV$, then let ``$\mu[\bb]$'' denote the measure
of the associated cylinder set:

  \[ \mu[\bb] \ := \ \mu\set{\ba \in \gA^\Zahl}{\ba_\sV \ = \ \bb } \]

  \item Suppose $\sV \subset \Zahl$ and  $k \in
\Zahl$ are such that  $(k + \sV) \ \subset \ \sU$.  If $\bb \ := \ \mtrx{b_v}{v \in \sV}{}$
is an element of $\gA^\sV$, then let $\bb'$ be the ``shift'' of $\bb$
by $k$: that is, $\bb' \ := \ \mtrx{b'_v}{v \in (k+\sV)}{}$, where,
for all $v \in \sV,  \ \ b'_v \ = \ b_{v-k}$.  Then define:

   \[ \mu_\sU[\bb] \ := \
     \mu_\sU\set{\ba \in \gA^\sU}{\ba_{(k + \sV)} \ = \ \bb'} \]

  (because $\mu_\sU$ is locally stationary, it doesn't matter which $k$ we use
in this definition, if more than one $k$ is available)

\eitem
  
  The {\bf Markov Extension} of $\mu_\sU$ is the probability
measure $\mu_{m\!r\!k}$, where, for any $N \geq U$, and $\bb \in \gA^\CC{0..N}$,

 \[ \mu_{m\!r\!k}[\bb] \ := \ \mu_\sU\lb[\bb_\sU\rb] \cdot
        \prod_{k=1}^{N-U} 
    \mu_\sU\lb[\frac{\bb_\CC{k...U\!+\!k}}{ \bb_\CO{k...U\!+\!k} }\rb] \]
  
  Here, $\CO{k...U\!+\!k} \ := \ \{k,\ k\!+\!1,\ \ldots,\ k\!+\!U\!-\!1\}$, while
 $\CC{k...U\!+\!k} \ := \ \{k,\ k\!+\!1,\ \ldots,\ k\!+\!U\}$, and
$\mu_\sU\lb[\frac{\bb_\CC{k...U\!+\!k} }{ \bb_\CO{k...U\!+\!k} }\rb]$ is 
the {\bf conditional probability}:

 \[ \mu_\sU\lb[\frac{\bb_\CC{k...U+k}}{ \bb_\CO{k...U+k} }\rb]
 \ \ := \ \  \frac{ \mu_\sU\lb[\bb_\CC{k...U+k} \rb] }  
                  { \mu_\sU\lb[\bb_\CO{k...U+k} \rb] } \] 

 $\mu_{m\!r\!k}$ is a stationary probability measure on $\gA^\Natur$.
Define the probabilities of cylinder sets indexed by negative
coordinates by simply {\em shifting} them into the positive
domain.  Thus, $\mu_{m\!r\!k}$ is defined on all cylinder sets in $\gA^\Zahl$.
It is straightforward to check that the probability measure thus
defined is stationary, and that its marginal projection upon $\gA^\sU$
is equal to $\mu_\sU$.

\breath
   
   This construction indicates that a stationary extension of the
measure $\mu_\sU$ always exists.  In general, there may be many such
extensions.  Intuitively, $\mu_{m\!r\!k}$ is an extension built so as
to provide the maximum amount of ``random choice'' at each successive
coordinate.  Hence, the following variational principle is not too
surprising:

\subparagraph{Theorem:}  {\tt Maximal Entropy Property}

{\em
 Of all the different stationary extensions of $\mu_\sU$ that exist,
$\mu_{m\!r\!k}$ is the one possessing the largest {\bf process entropy},
which we define as:

  \[ H\lb(\mu_{m\!r\!k}\rb) \ \ := \ \ 
     \lim_{N \goto \oo} \frac{-1}{N} \sum_{\ba \in \gA^\CC{1...N}} 
               \mu_{m\!r\!k}[\ba] \log_2 \lb(\mu_{m\!r\!k}[\ba]\rb) \]
}
\Proof \ \ \  See, for example, \cite{Schlijper}. \qed

\breath

 Under what circumstances do {\em ergodic} extensions of $\mu_\sU$
exist?  Can we build an extension measure which is supported only on
{\em periodic} words of some fixed periodicity?  Also, what happens if
$\sU$ is {\em not} just an interval inside $\Zahl$?

\subsection{Extension on Lattices}

  Now, let $D > 0$, and let $\Zahl^D$ be a $D-$dimensional lattice.
Then $\gA^{\Zahl^D}$ is the space of
{\em $D-$dimensional configurations} on $\gA$.  If $\bk \in \Zahl^D$,
then the {\bf shift by $\bk$} is the map $\shift{\bk}:\gA^{\Zahl^D} \into 
\gA^{\Zahl^D}$ so that, if $\ba := \mtrx{a_\bn}{\bn \in \Zahl^D}{}$,
then $\shift{\bk} \ba \ \  := \ \mtrx{a'_\bn}{\bn \in \Zahl^D}{}$,
where $a'_\bn \ = \ a_{\bn-\bk}, \ \ \forall \bn \in \Zahl^D$.
\label{shift.defn}

 A {\bf stationary stochastic
process} is a probability measure $\mu$ on $\gA^{\Zahl^D}$ that
is {\bf invariant} under all shift maps.  That is, if
$\sV \subset \Zahl^D$ is any finite subset, and $\bb \in \gA^\sV$,
then for any $\bk \in \Zahl^D$,

 \[ \mu\lb[\shift{\bk}(\bb) \rb] \ \ = \ \ \mu[\bb] \]

  If $\sU \subset \Zahl^D$, and $\bk \in \Zahl^D$,
then define $\shift{\bk}\sU \ = \ \sU+ \bk$, and
define $\shift{\bk}:\gA^{\sU} \into 
\gA^{\sU+\bk}$ so that, if $\ba := \mtrx{a_\bn}{\bn \in \sU}{}$,
then $\shift{\bk} \ba \ \  := \ \mtrx{a'_\bn}{\bn \in \sU+\bk}{}$,
where $a'_\bn \ = \ a_{\bn-\bk}, \ \ \forall \bn \in \sU+\bk$.
A probability measure $\mu_\sU$ on $\gA^\sU$
is {\bf locally stationary} if  for any $\sV \subset \sU$, any
$\bb \in \gA^\sV$, and any $\bk \in \Zahl^D$ so that $\shift{\bk}{\sV} \subset \sU$
also,

  \[ \mu_\sU\lb[\shift{\bk}(\bb) \rb] \ \ = \ \ \mu_\sU[\bb] \]

\subparagraph{The Extension Problem:}
{\em Given a locally stationary
measure $\mu_\sU$ upon $\gA^\sU$, can we {\bf extend} it to a stationary
stochastic process $\mu$ on $\gA^{\Zahl^D}$, so that $\marg{\sU}{\mu} \ \ = \ \ \mu_\sU$?}

\breath

  The Extension Problem does not always have solutions, as examples in
Section \ref{counter} will show.  If we {\em can} solve the Extension
Problem, can we construct an extension which is ergodic?  (quasi)
Periodic or (weakly) mixing?

\subsection{Extension on Group Modules}

  Now, let $\dG$ be an arbitrary group, and let $\sM$ be a {\bf
$\dG$-module}: \  an arbitrary set equipped with a $\dG-$action.
A few examples of this to keep in mind:

\bitem
  \item $\sM \ := \ \Zahl^D$ and $\dG \ := \ \Zahl^D$, also, acting upon $\sM$
by {\bf translation}.

  \item $\sM \ := \ (\Zahl/P_1) \dirsum (\Zahl/P_2) \dirsum \ldots
(\Zahl/P_D)$, and $\dG \ := \ \Zahl^D$ acts upon $\sM$ by {\bf translation}
with periodic boundary conditions.

  \item $\dG$ is an arbitrary group, $\dH$ an arbitrary subgroup, and
$\sM \ := \ \dG/\dH$ is the set of {\bf right cosets}.  $\dG$ acts upon
$\sM$ by multiplication:  if $\bg \in \dG$ and $(\bk\dH) \in \sM$, then
$\bg.(\bk\dH) \ := \ (\bg.\bk)\dH$.  (Every {\bf transitive}
$\dG-$module is of this type, and {\em every} $\dG-$module can be written
as a disjoint union of transitive $\dG$-modules.)
\eitem

Let $\gA^\sM$ be the space of
{\em $\sM-$indexed configurations} on $\gA$.  If $\bg \in \dG$
then the {\bf shift by $\bg$} is the map $\shift{\bg}:\gA^{\sM} \into 
\gA^{\sM}$ so that, if $\ba := \mtrx{a_m}{m \in \sM}{}$,
then $\shift{\bg} \ba \ \  := \ \mtrx{a'_{m}}{m \in \sM}{}$,
where $a'_m \ = \ a_{\bg^{-1}.m}, \ \forall m \in \sM$.

 A {\bf $\dG$-invariant stochastic
process} is a probability measure $\mu$ on $\gA^{\sM}$ that
is {\bf invariant} under the shift action of $\dG$.  That is, if
$\sV \subset \sM$ is any finite subset, and $\bb \in \gA^\sV$,
then for any $\bg \in \dG$,

 \[ \mu\lb[\shift{\bg}(\bb) \rb] \ \ = \ \ \mu[\bb] \]

  If $\sU \subset \sM$ and $\bg \in \dG$, then define
$\shift{\bg}\sU \ = \ \bg.\sU \ = \ \set{\bg.u}{u\in\sU}$, 
and define $\shift{\bg}:\gA^{\sU} \into 
\gA^{\bg.\sU}$ so that, if $\ba := \mtrx{a_u}{u \in \sU}{}$,
then $\shift{\bg} \ba \ \  := \ \mtrx{a'_{u}}{u \in \bg.\sU}{}$,
where $a'_u \ = \ a_{\bg^{-1}.u}, \ \forall u \in \bg.\sU$.
A probability measure $\mu_\sU$ on $\gA^\sU$
is {\bf locally stationary} if  for any $\sV$ subset $\sU$, any
$\bb \in \gA^\sV$, and any $\bg \in \dG$ so that $\shift{\bg}{\sV} \subset \sU$
also,

  \[ \mu_\sU\lb[\shift{\bg}(\bb) \rb] \ \ = \ \ \mu_\sU[\bb] \]
               
  Again, we ask:

\subparagraph{The (group module) Extension Problem:}
{\em  Given a locally stationary
measure $\mu_\sU$ upon $\gA^\sU$, can we {\bf extend} it to a stationary
stochastic process $\mu$ on $\gA^{\sM}$, so that $\marg{\sU}{\mu} \ \ = \ \ \mu_\sU$?}

\breath

  If $\sM = \Zahl^D = \dG$, then this is just the Extension Problem on a
$D-$dimensional lattice.  If $\sM :=  (\Zahl/P_1) \dirsum (\Zahl/P_2)
\dirsum \ldots (\Zahl/P_D)$ and $\dG := \Zahl^D$, then a $\dG-$invariant
measure on $\gA^\sM$ is ``equivalent'' to a stationary stochastic
process on $\gA^{\Zahl^D}$ which is supported only on {\em periodic}
configurations with fundamental domain $\CO{0...P_1} \x 
\CO{0...P_2} \x \ldots \x \CO{0...P_D}$.  In Section \ref{envelope}, we will
demonstrate that, if $\sU \subset
\CO{0...P_1} \x 
\CO{0...P_2} \x \ldots \x \CO{0...P_D} \subset \Zahl^D$ is some
``small enough'' domain, then any locally stationary measure
$\mu_\sU$ can be identified with a locally invariant measure $\mu_{\sU'}$,
where $\sU' \subset \sM$ is the obvious ``representation'' of $\sU$
inside $\sM$.

\subsection{Organization of this paper}

  In \S\ref{sect.appl}, we motivate the Extension Problem by
discussing applications to the {\bf Invariant Measure Problem} for
subshifts of finite type and cellular automata.  In \S
\ref{counter}, we show that the Extension Problem is not trivial
by providing examples of locally stationary measures which {\em
cannot} be extended.  These examples imply two necessary
conditions for extendibility: the {\bf Entropy Condition} and the {\bf
Tiling Condition}.

  In \S\ref{sect.harm}, we review basic harmonic analysis on
configuration space, treating it as a compact abelian group, and
characterise the Extension Problem in terms of
constructing a suitable set of Fourier coefficients.  We use this
in \S\ref{finite}, where we consider extension on {\em finite}
$\dG$-modules, and show that, if $\nu$ is an extendible measure with
full support, and $\mu$ is ``close enough'' to $\nu$, then $\mu$ is
also extendible.  A similar result can be developed for
constructing {\em periodic} extensions, but first we need a
tool to ``reduce'' the Extension Problem on an infinite
module to an extension problem on a suitably chosen finite module,
which we develop in \S\ref{envelope}, via the concept of ``envelopes''.

  In \S\ref{embed}, we show that an extendible, locally
stationary measure with full support can be ``embedded'' in any
ergodic $\Latt$-dynamical system, in the sense that it is a marginal
projection of a stationary $\Latt$-process generated by a partition on
that system.  

  In \S\ref{sect.per}, we combine the results of \S\ref{finite} and 
\S\ref{envelope} to investigate when a measure has an 
almost-surely periodic extension, and provide examples of measures
which {\em never} have periodic extensions, as well as measures which
{\em only} have periodic extensions.  Then we use the results of
\S\ref{embed} to show that ``almost all'' extendible measures
have extensions which are ergodic, mixing, weakly mixing, or
quasiperiodic.

  In \S \ref{decide}, we show that the Extension Problem is,
in general, formally undecidable.

\subsection{Preliminaries and Notation}

  If we treat $\gA$ as a {\bf discrete} topological space, and endow
$\gA^\sM$ with the Tychonoff product topology, then $\gA^\sM$ is a
compact, metrizable space.  If $\sM$ is finite, then $\gA^\sM$ is
finite and discrete.  If $\sM$ is infinite, then $\gA^\sM$ is
uncountable and {\bf totally disconnected}.

   The topology on $\gA^\sM$ is generated by 
{\bf cylinder sets}.  If $\sU \subset \sM$ is finite, and $\bb \in \gA^\sU$,
then the associated cylinder set is:

   \[ \set{\ba \in \gA^\sM}{\ba_\sU \ = \ \bb } \]

  Here, by ``$\ba_\sU$'' we mean the element $\mtrx{a_u}{u \in \sU}{}$, 
where $\ba =  \mtrx{a_m}{m \in \sM}{}$.  Normally, we will use the
symbol ``$\bb$'' to denote both the word $\bb$ and the cylinder set
it induces ---the distinction will be clear from context.  For example, 
if $\mu$ is some measure, then ``$\mu[\bb]$'' indicates the measure
of the cylinder set defined by $\bb$.

 Whenever we speak of measures on $\gA^\sM$, we will mean measures
on the {\bf Borel sigma-algebra} generated by the product topology.

\breath

  If $\sM$ is a $\dG-$module, then $\invMeas{\dG}{\gA^\sM}$ is the space of all
$\dG-$invariant probability measures on $\gA^\sM$.  This is a convex
subset of $\Meas{\gA^\sM}$, the space of {\em all} {\bf probability
measures} on $\gA^\sM$, which, in turn, is a convex subset of
the real vector space $\Meas{\gA^\sM; \ \Real}$ of {\bf real-valued
measures} on $\gA^\sM$.  

  The elements of $\Meas{\gA^\sM; \ \Cplx}$ ({\bf complex-valued
measures} on $\gA^\sM$) act as linear functionals on $\bC(\gA^\sM; \ \Cplx)$
(the Banach space of {\bf complex-valued, continuous functions}).
This induces a {\bf weak$-*$ topology} on 
$\Meas{\gA^\sM; \ \Cplx}$, making it into a locally
convex topological vector space. 

  $\invMeas{\dG}{\gA^\sM}$ is a compact subset of $\Meas{\gA^\sM; \ \Cplx}$
under this topology.

  When $\dG = \sM = \Latt$, we will refer to $\invMeas{\dG}{\gA^\sM}$ as
``$\statMeas{\gA^\Latt}$''.

\breath
  
  If $\sU \subset \sM$, then $\invMeas{\dG}{\gA^\sU}$ is the space of all
locally $\dG-$invariant probability measures on $\gA^\sU$.
$\extMeas{\gA^\sU}$ is the set of all {\bf extendable probability measures}:
measures which can be extended to a $\dG-$invariant measure on $\gA^\sM$.
Notice that:

\begin{quote}
  {\em $\extMeas{\gA^\sU}$ is a compact, convex subset of 
$\Meas{\gA^\sU; \Cplx}$.}
\end{quote}

  This is because 
the marginal projection map $\pr{\sU}^*:\Meas{\gA^\sM; \ \Cplx} \into
\Meas{\gA^\sU; \ \Cplx}$ is {\bf linear} and {\bf continuous},
and $\extMeas{\gA^\sU}$ is simply the image of the compact,
convex subset $\invMeas{\dG}{\gA^\sM}$ under $\pr{\sU}^*$.

\section{Applications \label{sect.appl}}

\subsection{Subshifts of Finite Type \label{finite.type}}

 Let $\sU \subset \Latt$ be finite, and suppose that $\gW \subset \gA^\sU$ is some
set of ``admissible'' $\sU$-words.   The {\bf subshift of finite type}
defined by $\gW$ is the closed, shift-invariant subset of
$\gA^\Latt$:

  \[ \inn{\gW} \ := \ \set{\ba \in \gA^\Latt}
       { \forall \bn \in \Latt, \ \ \ba_{\sU+\bn} \in \gW} \]

  One-dimensional subshifts of finite type were first studied by Parry
\cite{ParrySoFT} and Smale \cite{Smale}; \ excellent recent
introductions are \cite{LindMarcus} and \cite{Kitchens}.  Higher
dimensional subshifts are closely related to tilings \cite{Mozes1},
\cite{Mozes2},\cite{Radin}, and involve many additional subtleties; \ see,
for example \cite{MarkleyPaul1},\cite{MarkleyPaul2}.  Of particular
interest is

\subparagraph{The Nontriviality Problem:}
{\em For a given set $\gW$, is the corresponding
set $\inn{\gW}$ is even nonempty? }

\breath

   The Nontriviality Problem is known to be {\bf formally
undecidable}; see \cite{RMR}, \cite{BerDomino}, or \cite{KitS}.

\Theorem{}
{ Let $\sU$ and $\gW$ be as above.  $\inn{\gW}$ is nontrivial if and
only if there is some locally stationary probability measure $\mu_\sU$
on $\gA^\sU$, with $\supp{\mu_\sU} \subset \inn{\gW}$, such that
$\mu_\sU$ has a stationary extension.}

\begin{thmproof}
  Suppose that such a $\mu_\sU$ existed, and let $\mu$ be a stationary extension.
Clearly, any $\mu-$generic configuration in $\gA^\Latt$ must satisfy the
membership criteria of $\inn{\gW}$.  Hence, $\inn{\gW}$ must be nonempty.

  Conversely, if $\inn{\gW}$ was nonempty, then by the Krylov-Bogoliov theorem
\cite{Walters}, there are stationary probability measures whose 
support is contained in $\inn{\gW}$.  Let $\mu$ be one of these measures, and
let $\mu_\sU := \marg{\sU} \mu$.  Then $\supp{\mu} \subset \gW$.
\end{thmproof}

  Let $\extMeas{\gW}$ be the set of extendible measures supported
on $\gW$.

\Corollary{}
{  It is formally undecidable whether, for a given subset $\gW \subset \gA^\sU$,
the set $\extMeas{\gW}$ is nonempty.}
\qed

\breath

  However, it is easily decidable whether $\statMeas{\gW}$ itself is
nonempty.  The set of all real-valued measures supported on $\gW$ is
a finite-dimensional vector space, and the stipulation that an element
of this vector space be a locally stationary probability measure takes
the form of a finite system of linear equations and inequalities; \
solving such a system is a decidable problem.

\subsection{Cellular Automata \label{CA} }

 Let $\sU \subset \Latt$ be finite (metaphorically speaking, $\sU$ is
a ``neighbourhood of zero'') and let $\phi:\gA^\sU \into \gA$.  For
every $\bn \in \Latt$, define $\phi_\bn := \phi \circ \shift{-\bn} :
\gA^{\sU+\bn} \into \gA$.

 The {\bf cellular automata} determined by $\phi$ is then the function 
$\Phi:\gA^\Latt  \into 	\gA^\Latt$ sending 
$\mtrx{a_\bn}{\bn \in \Latt}{} \mapsto 
\mtrx{\phi_\bn \lb(\ba_{\sU+\bn} \rb)}{\bn \in \Latt}{}$.
$\phi$ is called the {\bf local transformation rule} for $\Phi$.
Cellular automata were first investigated by Von Neumann \cite{vonNeumannCA}
and Ulam \cite{UlamCA}, and later extensively studied by Hedlund
\cite{HedlundCA}, Wolfram \cite{WolframBook}, and others; \ more
recent surveys are \cite{PhysicaCA1},\cite{PhysicaCA2},\cite{DelormeMazoyer},
\cite{GolesMartinez}.

\breath

  Any cellular automaton on $\Latt$ can be represented by a subshift of finite
type on $\Latt \x \Zahl$.  Simply define

  \[ \tlsU \ := \ \lb( \sU \x \{0\} \rb) \disj
   \lb\{ ( \underbrace{0,0,\ldots,0}_D , \ 1 ) \rb\} \]

  and then set $\tlgW \ := \ \set{\ba \in \gA^\tlsU}{
    a_{( 0,0,\ldots,0, \ 1 )} \ = \ 
   \phi \lb( \ba_{\lb( \sU \x \{0\} \rb)} \rb)}$

  If $\ba \in \gA^{\Latt \x \Zahl}$, then $\ba$ can be seen as a $\Zahl-$indexed
sequence of configurations in $\gA^\Latt$.  Clearly, $\ba$ is in $\inn{\tlgW}$ if and
only if this sequence describes the $\Phi-$orbit of some point in $\gA^\Latt$.

  Of course, unless $\Phi$ is surjective, not every element of $\gA^\Latt$ will
necessarily have a $\Phi-$preimage, and thus, not every element can
appear in such a $\Zahl-$indexed sequence of configurations.  We can
obviate this difficulty by concentrating on the {\bf center} of
the dynamical system $(\gA^\Latt, \ \Phi)$.

  If $X$ is any compact space, and $T:X \into X$ continuous, then
the {\bf nonwandering set}, $\Omega(X,T)$ is the set of
all points $x \in X$ which are {\bf regionally recurrent}:  for
any neighbourhood $U$ of $x$, there is some $n \in \Natur$ so that
$T^n(U) \intsct U \not= \emptyset. \ \ \ \Omega(X,T)$ is a compact
$T-$invariant subset, so we can look at the restricted dynamical system
$\lb(\Omega(X,T), \ T_{|\Omega(X,T)} \rb)$ ---however, not all
elements of $\Omega(X,T)$ will be regionally recurrent under
$T$, when seen in the subspace topology (see \cite{Walters} for
an example) ---hence, $\Omega^2(X,T) := 
\Omega\lb(\Omega(X,T), \ T_{|\Omega(X,T)} \rb)$ may be a proper
subset. 

  By transfinite induction,  for any countable ordinal number $\alpha$,
define \ $\Omega^{\alp+1}(X,T) := 
\Omega\lb(\Omega(X,T), \ T_{|\Omega^\alp(X,T)} \rb)$, 
and, if $\gamma$ is a limit ordinal, define $\Omega^{\gam}(X,T)$\ $:=
 \Intsct_{\alp < \gam} \Omega^{\alp}(X,T)$.  Since $X$ is compact,
 this descending sequence of compact subsets must become constant at
 some countable ordinal $\alp$, so that $\Omega^{\alp+1}(X,T)
 =\Omega^{\alp}(X,T)$.  The {\bf center} of $(X,T)$, defined $\bZ(X,T)
 := \Omega^{\alp}(X,T)$, is nonempty, compact, and $T-$invariant.  If
 $\mu$ is any $T-$invariant Radon measure on $X$, then $\supp{\mu}
 \subset \bZ(X,T)$.

  So, treat $(\gA^\Latt, \ \Phi)$ as a compact topological dynamical
system, and let $\bZ(\Phi)$ be its center. 
The restricted
map $\Phi_|:\bZ(\Phi) \into \bZ(\Phi)$ is surjective, so every element
in $\bZ(\Phi)$ appears in some $\Zahl-$indexed sequence of $\gA^\Latt-$configurations
admissable to $\tlgW$.

\subparagraph{The Invariant Measure Problem:}
{\em  Given a local transformation rule $\phi:\gA^\sU \into \gA$,
describe the set of $\Phi-$invariant, stationary measures on $\gA^\Latt$.}

\breath

  Suppose that we represent the cellular automata as a subshift of
finite type in the aforementioned way, and suppose that $\mu_\tlsU$ is a
locally stationary probability measure on $\gA^\tlsU$.  It is easy to
verify that a stationary extension of $\mu_\tlsU$ to $\gA^{\Latt \x \Zahl}$
is equivilant to a $\Phi-$invariant, stationary measure on $\gA^\Latt$.

\section{Caveats and Counterexamples \label{counter}}

\subsection{Nonextendability in $\Zahl$;  The Entropy Metric} 

  The following counterexample, which first appeared in \cite{DJRW},
shows that, even in $\Zahl$, locally stationary measures are not
always extendible, when the initial domain is ``disconnected''.

  Suppose that $\sU := \{0,1,3\}$.  If $\mu_\sU$ is a probability
measure on $\gA^\sU$, then we can treat the functions $\pr{0}, 
\pr{1}$, and $\pr{3}$ as {\em random variables} ranging over the domain
$\gA$.  So, let $\mu_\sU$ be any probability measure
on $\gA^\sU$ such that:

\bitem
  \item (A) $\pr{0} = \pr{1}, \  \mu_\sU-$almost-surely.
  \item (B) $\pr{0}$ and $\pr{3}$ are {\em independent} as random variables.
(thus $\pr{1}$ and $\pr{3}$ are also independent.)
\eitem

  To ensure $\mu_\sU$ is locally stationary, it suffices to require
only that the random variables $\pr{0}$, $\pr{1}$, and $\pr{3}$ are
{\em identically distributed}.  

  The measure $\mu_\sU$ cannot be extended even to a locally
stationary measure on $\gA^\CC{0..3}$, much less a stationary measure
on $\gA^\Zahl$.  To see this, suppose that $\mu_\CC{0..3}$ was a
locally stationary extension.  Then condition (A) defining $\mu_\sU$
implies that, as random variables on the probability space
$\lb(\gA^\CC{0..3}, \ \mu_\CC{0..3} \rb), \ \ \pr{0} = \pr{1} \ =
\pr{2} \ = \ \pr{3}$.  But by condition (B),
$\pr{0}$ and $\pr{3}$ are independent ---a contradiction.

\breath

  This example can be understood as part of a more general phenomenon.
If $\sS$ is any set, and $\mu$ is any probability measure on
$\gA^\sS$, then $\mu$ induces an {\bf entropy metric}, $D_\mu$, on the
set $\Fin{\sS}$ of all finite subsets of $\sS$.  If $\sU, \sV \subset
\sS$ are finite, then define

  \[ H_\mu [\sU | \sV] \ := 
  \ - \sum_{\bb \in \gA^\sV} \sum_{\ba \in \gA^\sU} \mu[\ba|\bb] \log_2(\mu[\ba|\bb]), \]

 \[ \btxt{ where \etxt} \ \  \mu[\ba|\bb] \ := \ 
  \frac{\mu\set{\bc \in \gA^\sS}{\bc_\sU = \ba \btxt{ \ and \ \etxt} 
              \bc_\sV = \bb }}
                                 {\mu\set{\bc \in \gA^\sS}{\bc_\sV = \bb}}. \]

\[ \btxt{  Then define: \etxt} \ \ 
  D_\mu[\sU, \ \sV] \ := \  H_\mu [\sU | \sV] +  H_\mu [\sV | \sU]. \]

  It is easy to check that $D_\mu$ is a {\em metric} on $\Fin{\sS}$.  Furthermore,
if $\sS$ is a $\dG-$module, and $\mu$ is a $\dG-$invariant measure, then
$D_\mu$ is a $\dG-$action invariant metric.  If $\sS$ is a {\em subset}
of some $\dG-$module, and $\mu$ is a {\em locally} $\dG-$invariant
measure, then $D_\mu$ is a ``locally'' $\dG$-invariant metric, in the obvious
sense.

\breath

  Now, suppose $\sM$ is a $\dG-$module, $\sU \subset \sM$, and $\mu_\sU$ is a
  locally $\dG-$invariant
measure on $\gA^\sU$.  If $\mu$ is to be an invariant extension of 
$\mu_\sU$, then it must satisfy the condition:

\bquote
{\em For every $\sV, \sW \in \Fin{\sU}$, and every $\bg \in \dG, \ \ 
 D_\mu[\bg.\sV, \ \bg.\sW] \ = \ D_{\mu_\sU} [\sV, \sW].$}
\equote

  Hence, $D_\mu$ is forced to take certain values on a subset of
$\Fin{\sM}$.  The question is: can we define $D_\mu$ in the {\em rest}
of $\Fin{\sM}$ so that it is a metric?  If we cannot, then it is
impossible to extend $\mu$.

\breath

  In the aforementioned counterexample, $D_{\mu_\sU} \lb[ \{0\},\  \{1\} \rb]
= 0$.  Thus, if $\mu$ was an extension of $\mu_\sU$, we would have:

  \[ D_\mu \lb[ \{0\},\  \{1\} \rb] \ = \
     D_\mu \lb[ \{1\},\  \{2\} \rb] \ = \
     D_\mu \lb[ \{2\},\  \{3\} \rb] \ = \ 0 \]

  and hence, $D_\mu \lb[ \{0\},\  \{3\} \rb] \ = \ 0$.  But we know
that  $D_\mu \lb[ \{0\},\  \{3\} \rb] \ > \ 0$, because $\pr{0}$ and
$\pr{3}$ are independent random variables.  Hence, no such
extension $\mu$ can exist.

\subsection{Nonextendability in $\Zahl^D$;  The Tiling Condition
\label{tiling.condition}} 

 In the previous counterexample, it seems the problem was that
the domain $\sU$ was not ``connected''.  However, in $\Zahl^2$,
extendability can fail even when $\sU$ is a $2 \x 2$ box.

 Suppose $\sU \subset \Latt$, and $\mu_\sU \in \statMeas{\gA^\sU}$.  
The {\bf support} of $\mu_\sU$ is some subset $\supp{\mu_\sU} \subset \gA^\sU$;
 \  let $\inn{\supp{\mu_\sU}}$ be the {\bf subshift of finite type} defined
by $\supp{\mu_\sU}$.  If
$\mu \in \statMeas{\gA^\Latt}$ is a stationary extension of $\mu_\sU$,
then any $\mu-$generic configuration $\ba \in \gA^\Latt$ must
be an element of $\inn{\supp{\mu_\sU}}$. 

Thus, we have:

\subparagraph{The Tiling Condition:}
  {\em  $\mu_\sU$ cannot be extendible unless $\inn{\supp{\mu_\sU}}$ is
nontrivial.}

\breath

  Intuitively, the configuration $\ba$ determines a
{\bf tiling} of $\Latt$ by elements in $\supp{\mu_\sU}$: \
 for any $\bk \in \Latt,
\ \ \ba_{(\bk+\sU)}$ is an element of $\supp{\mu_\sU}$.  

  For example, suppose that $D := 2, \ \ U := \CC{0..1} \x \CC{0..1}$,
and $\gA := \{0,\ 1,\ 2\}$.  Elements of $\gA^\sU$ are thus $2 \x 2$ words
in $\gA$. 

\begin{figure}[htbp]
\[\begin{array}{ccccccc}
\ddots	& \vdots & \vdots & \vdots & \vdots & \vdots & \\
\ldots	&0 &0 &2 &1 &1 & \ldots \\
\ldots	&1 &0 &1 &0 &1 & \ldots \\
\ldots	&2 &1 &2 &1 &2 & \ldots \\
\ldots	&0 &1 &1 &1 &1 & \ldots \\
\ldots	&0 &1 &0 &0 &1 & \ldots \\
	& \vdots & \vdots & \vdots & \vdots & \vdots & \ddots \\
\end{array}\]
\caption{A configuration of letters \label{fig1}}
\end{figure}

\begin{figure}[htbp]
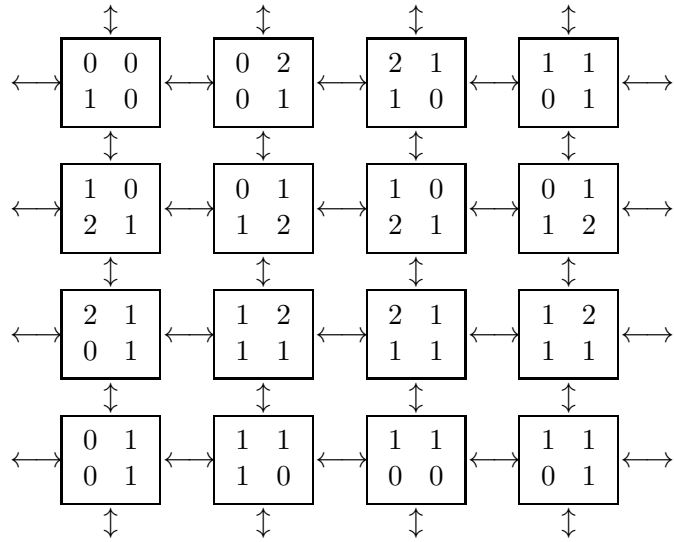

\[\begin{array}{cc cc cc ccc}     
 & \updownarrow & & \updownarrow & & \updownarrow &  & \updownarrow\\ 
\connect & \tilebox{cc}{ 0 & 0  \\ 1 & 0 } &    
\connect & \tilebox{cc}{ 0 & 2  \\ 0 & 1 } &    
\connect & \tilebox{cc}{ 2 & 1  \\ 1 & 0 } &    
\connect & \tilebox{cc}{ 1 & 1  \\ 0 & 1 } &    
\connect \\   
 & \updownarrow & & \updownarrow & & \updownarrow &  & \updownarrow\\ 
\connect & \tilebox{cc}{ 1 & 0  \\ 2 & 1 } &    
\connect & \tilebox{cc}{ 0 & 1  \\ 1 & 2 } &    
\connect & \tilebox{cc}{ 1 & 0  \\ 2 & 1 } &    
\connect & \tilebox{cc}{ 0 & 1  \\ 1 & 2 } &   
\connect \\   
 & \updownarrow & & \updownarrow & & \updownarrow &  & \updownarrow\\ 
\connect & \tilebox{cc}{ 2 & 1  \\ 0 & 1 } &    
\connect & \tilebox{cc}{ 1 & 2  \\ 1 & 1 } &    
\connect & \tilebox{cc}{ 2 & 1  \\ 1 & 1 } &    
\connect & \tilebox{cc}{ 1 & 2  \\ 1 & 1 } &   
\connect \\   
 & \updownarrow & & \updownarrow & & \updownarrow &  & \updownarrow\\ 
\connect & \tilebox{cc}{ 0 & 1  \\ 0 & 1 } &    
\connect & \tilebox{cc}{ 1 & 1  \\ 1 & 0 } &    
\connect & \tilebox{cc}{ 1 & 1  \\ 0 & 0 } &    
\connect & \tilebox{cc}{ 1 & 1  \\ 0 & 1 } &   
\connect \\   
 & \updownarrow & & \updownarrow & & \updownarrow &  & \updownarrow\\ 
\end{array}\] 
\caption{The corresponding assignment of matrices. \label{fig2}}
\end{figure}

  Choosing a configuration in $\gA^{\Zahl^2}$ is equivalent to
assigning a $2 \x 2$ matrix to each point in the lattice, so that
adjacent sides agree.  For example, the configuration in Figure \ref{fig1}
is equivalent to the assignment of Figure \ref{fig2}

\begin{figure}
\includegraphics[scale=2]{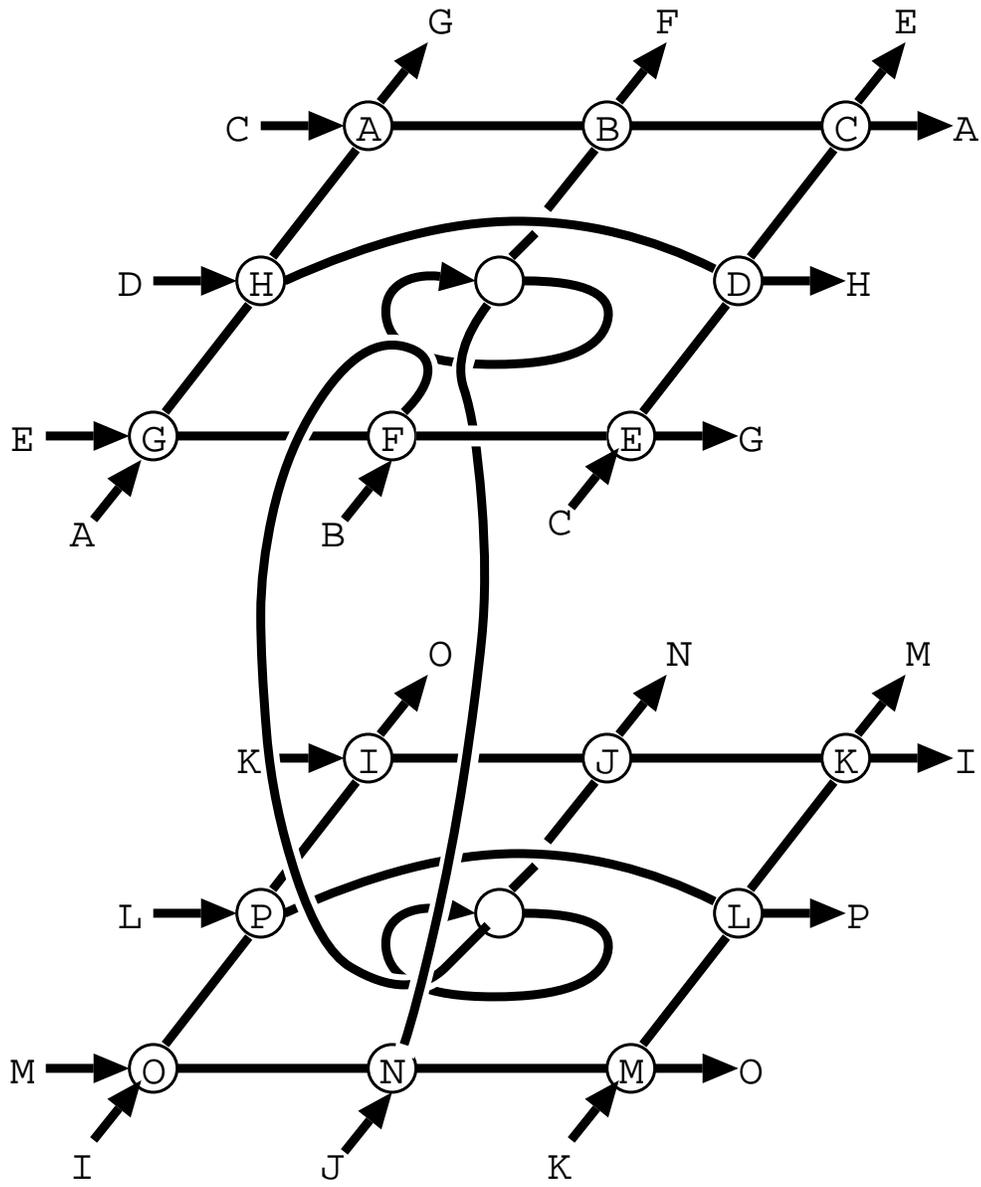}
\caption{A ``pseudolattice''.\label{Fig:pseudolattice}}
\end{figure}

   We will define a locally stationary measure $\mu_\sU$ so that
$\supp{\mu_\sU}$ cannot tile $\Zahl^2$ in this manner.  We will
do this by explicitly constructing $\supp{\mu_\sU}$ to tile
a {\em different} space instead ---a kind of ``pseudolattice''
(see Figure \ref{Fig:pseudolattice}).

   Stack two $3 \x 3$ grids on top of one another, and 
then ``break'' the connection between the central element of each
level, and its southern, eastern, and western neighbours.
Cross-connect the eastern and western neighbours with each other.
Connect the southern neighbour to the central element of the
level {\em above}, and we connect the central element of this level to the
southern element of the level {\em below}.  We also maintain the
connection between the central element and its northern neighbour,

\begin{figure}
\[\begin{array}{cc ccc cc}  
  & \updownarrow & & \updownarrow & & \updownarrow & \\
\connect & \tilebox{cc}{ 9 & 1 \\ 10 & 0 } & \connect & 
\tilebox{cc}{ 1 & 4 \\ 0 & 0 } & \connect & 
\tilebox{cc}{ 4 & 9 \\ 0 & 10 } & \connect \\  
 &  \updownarrow & & \updownarrow & & \updownarrow & \\ 
 \connect & \tilebox{cc}{ 10 & 0 \\ 11 & 7 } &  & 
\tilebox{cc}{ 0 & 0 \\ 6 & 6 } &  & 
\tilebox{cc}{ 0 & 10 \\ 7 & 11 } &  \connect \\  
 & \updownarrow & &              & & \updownarrow & \\ 
\connect & \tilebox{cc}{ 11 & 7 \\ 9 & 1 } & \connect & 
\tilebox{cc}{ 7 & 7 \\ 1 & 4 } & \connect & 
\tilebox{cc}{ 7 & 11 \\ 4 & 9 } &  \connect \\  
& \updownarrow & & \updownarrow & & \updownarrow & \\
\end{array}\]

\[\begin{array}{ccccccc}  
  & \updownarrow & & \updownarrow & & \updownarrow & \\
\connect & \tilebox{cc}{ 13 & 12 \\ 14 & 2 } & \connect & 
\tilebox{cc}{ 12 & 5 \\ 2 & 2 } & \connect & 
\tilebox{cc}{ 5 & 13 \\ 2 & 14 } & \connect \\  
 &  \updownarrow & & \updownarrow & & \updownarrow & \\ 
 \connect & \tilebox{cc}{ 14 & 2 \\ 15 & 6 } &  & 
\tilebox{cc}{ 2 & 2 \\ 7 & 7 } &  & 
\tilebox{cc}{ 2 & 14 \\ 6 & 15 } &  \connect \\  
 & \updownarrow & &              & & \updownarrow & \\ 
\connect & \tilebox{cc}{ 15 & 6 \\ 13 & 12 } & \connect & 
\tilebox{cc}{ 6 & 6 \\ 12 & 5 } & \connect & 
\tilebox{cc}{ 6 & 15 \\ 5 & 13 } &  \connect \\  
& \updownarrow & & \updownarrow & & \updownarrow & \\ 
\end{array}\] 
\caption{A configuration on the pseudolattice \label{fig4}}
\end{figure}

  Now we'll form a locally stationary measure which tiles {\em this}
space instead.  Consider the tiling portrayed in Figure \ref{fig4}.
Count every element of $\gA^{2 \x 2}$ as many times as it
appears in these two pictures.  There are $18$ tiles, and each one
appears exactly once.  Thus, each of the tiles shown gets a
probability of $\frac{1}{18}$.

  To show that $\mu_\sU$ is locally stationary, it suffices to check
that the {\em left columns} and {\em right columns}
have the same probability distribution, and that the {\em top} and {\em bottom
rows} have the same probability distribution.  This is easy to confirm.

  We claim that one simply {\em cannot} tile $\Zahl^2$ with
this collection of blocks.  For example, as soon as one lays down a
tile of the form $\minitilebox{0}{0}{6}{6}$, one is forced to place a
tile $\minitilebox{ 1 }{ 4 }{ 0 }{ 0 }$ immediately above it, since
this is the only tile which will ``match''.  Once one has done this,
one {\em must} place the tile $\minitilebox{ 9 }{ 1 }{ 10 }{ 0 }$ to the
left of $\minitilebox{ 1 }{ 4 }{ 0 }{ 0 }$, and the tile
$\minitilebox{ 4 }{ 9 }{ 0 }{ 10 }$ to its right.  So far, all the
tiles are compatible.  However, now, what tile shall we lay down below
$\minitilebox{ 9 }{ 1 }{ 10 }{ 0 }$?  To be compatible with
$\minitilebox{ 9 }{ 1 }{ 10 }{ 0 }$, this tile's top row should read
$\minitilebox{10 }{ 0 }{}{}$.  However, to be compatible with the tile
$\minitilebox{ 0 }{ 0 }{ 6 }{ 6 }$ to its immediate right, the tile's
right-hand side should read $\minitilebox{\ \ }{ 0 }{\ \ }{ 6}$.  There is no
tile in our collection which meets these two criteria.

\breath

  The Tiling Condition is necessary, but not sufficient.  To see this,
recall that the set $\extMeas{\sU}$ is {\bf closed} as a subset of
$\statMeas{\sU}$.  Thus, its complement is {\bf open}.  Hence, every
nonextendible measure is surrounded by a neighbourhood of
nonextendible measures.

 If $\bet$ is the {\em equidistributed} measure (assigning equal
probability to {\em every} element of $\gA^\sU$), and $\eps > 0$ is
small, then consider the measure:

    \[ \mu_\eps \ \ \ := \ \ (1 - \eps)\cdot \mu_\sU \ + \ \eps \cdot \bet \]

  $\mu_\eps$ is a convex combination of $\mu_\sU$ and $\bet$.
Since  $\eps > 0$, the support of $\mu_\eps$ is {\em all} of $\gA^\sU$.
Thus, $\mu_\eps$ always satisfies the Tiling Condition.
However, if $\eps$ is ``sufficiently small'', the measure $\mu_\eps $
will be inside the neighbourhood of nonextendible measures
around $\mu_\sU$.

\section{Harmonic Analysis of Extensions \label{sect.harm}}

\subsection{Configuration Space as a Compact Group}

  Solving the Extension Problem requires a good way of describing
measures, and Harmonic Analysis provides one.
To employ this approach, we must reconceive the configuration
space as a {\bf compact abelian topological group}.  Hence,
from now on, we will operate under the assumption that:

\bquote
{\em The alphabet $\gA$ is a finite abelian group.}
\equote

  The choice of group structure on $\gA$ is unimportant ---if $\gA$ has
$A$ elements, then the simplest choice is to let $\gA := \Zahl/A$.

  If we endow $\gA^\sM$ with the product group structure,
it is a compact abelian topological group.  What is its dual group?

 Let $\hgA$ be the dual group of $\gA$.  If $\sV \subset \sM$ is finite,
and, for all $v \in \sV, \ \ \chi_v \in \hgA$, then $\chi_v \circ \pr{v}: \gA^\sM \into \Torus{1}$ is the map taking the configuration $\mtrx{a_m}{m \in \sM}{}$
to the value $\chi_v(a_v)$.  (Here ``$\Torus{1}$'' is the {\bf unit circle} 
group.)

 \[ \btxt{ We will use the notation \ ``\etxt} \Tensor_{v \in \sV} \chi_v 
 \btxt{'' to refer to the map: \ \  \etxt}  \]

\beq
  \lb( \prod_{v \in \sV} \chi_v \circ \pr{v} \rb) : \gA^\sM & \into & \Torus{1} \\
           \mtrx{a_m}{m \in \sM}{} & \mapsto & 
                         \prod_{v \in \sV} \chi_v(a_v) \\
\eeq

  It is easy to verify the next theorem:

\Theorem{}
{ Let $\sM$ be any set.  The dual group of $\gA^\sM$ is the set: 

  \[ \set{ \Tensor_{v \in \sV} \chi_v }
 { \sV \subset \sM \btxt{ \  is any finite subset, and, for all \  \etxt} v \in \sV,   \ \chi_v \in \hgA. } \]
}
\qed

\subsection{The Fourier Transform}

  Now, if $\mu$ is a measure on $\gA^\sM$ , and $\chi \in \dual{\sM}$, then
the {\bf Fourier Coefficient} of $\mu$ at $\chi$ is defined:

  \[ \hmu_\chi \ \ = \ \ 
\inn{\mu,\ \chi} \ \ := \ \ \int_{\gA^\sM} \barchi \ d\mu \]

  The {\bf Fourier Transform} of $\mu$ is the function: \ 
$  \hmu: \dual{\sM}  \into  \Cplx$ so that
$\hmu_\chi \ \ = \ \  \inn{\mu,\ \chi}$.

 If $\Meas{\gA^\sM; \ \Cplx}$ is endowed with the {\bf total variation norm},
and $\bC( \dual{\sM}; \ \Cplx)$ is endowed with the {\bf uniform norm},
then the map

\beq
   \Four:  \Meas{\gA^\sM; \ \Cplx} & \into &  \bC( \dual{\sM}; \ \Cplx) \\
             \mu & \mapsto & \hmu \\
\eeq

  is an {\bf injective, bounded linear function} of norm 1 
\cite{Katznelson}.  Thus, the Fourier transform of $\mu$
totally characterizes it: if $\mu$ and $\nu$ are two measures, and
$\hmu = \hnu$, then $\mu = \nu$.

\subsection{Fourier Theory and (local) Stationarity}

  The shift action of $\dG$ upon $\gA^\sM$ induces a right action of
$\dG$ upon $\dual{\sM}$.  If $\bg \in \dG$, and $\chi \in \dual{\sM}$,
then define:

\beqn
\label{char.action}
 \chi.\bg \ = \ \chi \circ \shift{\bg^{-1}} 
\eeqn

 \[ \btxt{ Note that, if \ \ \etxt}
  \chi \ = \ \prod_{v \in \sV} \lb( \chi_v \circ \pr{v} \rb), \btxt{ \ \ then
 \ \  \etxt}
\chi.\bg \ = \ \prod_{v \in \sV} \lb( \chi_v \circ \pr{\bg.v} \rb) \]

\breath

  If $\sU \subset \sM$ is not closed under the $\dG-$action, then
there is no ``shift action'' on $\gA^\sU$.  However, we can still treat
$\dG$ as ``acting'' upon $\dual{\sU}$ in a certain limited capacity, as follows:

  Suppose $\sV \subset \sU$, and $\chi \ = \ \prod_{v \in \sV}\  \lb(\chi_v
\circ \pr{v}\rb)$.  Suppose that $\bg \in \dG$ is
such that $\bg.\sV \subset \sU$ also.  Then $\chi.\bg \ = \ \prod_{v
 \in \sV} \lb( \chi_v \circ \pr{\bg.v} \rb)$ is still an element of
 $\dual{\sU}$.
 
\Theorem{\label{stat.four}}
{
\blist
  \item  If $\mu \in \Meas{\gA^\sM}$, then $\mu$ is $\dG-$invariant if and only if,
for every $\chi \in \dual{\sM}$ and every $\bg \in \dG$, \ \ 
$ \inn{\mu, \ \chi} \ \ = \ \ \inn{\mu, \ \chi.\bg}$

  \item  If $\sU \subset \sM$, and $\mu \in \Meas{\gA^\sU}$, then $\mu$ is 
{\em locally} $\dG-$invariant if and only if,
for every $\chi \in \dual{\sU}$ and every $\bg \in \dG$ so that $\chi.\bg$ is
also in $\dual{\sU}$,  \ \ $\inn{\mu, \ \chi} \ \ = \ \ \inn{\mu, \ \chi.\bg}$. 
\elist
}
\begin{thmproof}
  We will prove \thmpart{2}, since \thmpart{1} clearly follows. 

\proofof{``$\implies$''} 
Let $\chi = \Tensor_{v \in \sV} \chi_v$, for some $\sV \subset \sU$.
Then a simple computation reveals:

 \[ \inn{\mu, \chi} \ \ = 
  \ \ \sum_{\ba \in \gA^\sV} \mu[\ba] \cdot 
      \barchi(\ba) \]

 Where, by ``$\mu[\ba]$'', we mean 
$\mu\set{\bb \in \gA^\sU}{\bb_{\sV} \ = \ \ba}$. Thus,

\beq
 \inn{\mu, \chi.\bg} 
 & =_{(1)} &
  \sum_{\ba \in \gA^{\bg.\sV}}
        \mu[\ba] \cdot \lb(\barchi\circ\shift{\bg^{-1}} (\ba)\rb) \\
 & =_{(2)} &
  \sum_{\ba \in \gA^{\sV}}
        \mu\lb[\shift{\bg}\ba\rb] \cdot \lb(\barchi\circ
	\shift{\bg^{-1}}\circ\shift{\bg} (\ba) \rb) \\
 & =  &
  \sum_{\ba \in \gA^{\sV}}
        \mu\lb[\shift{\bg}\ba\rb] \cdot \barchi (\ba)  \\
 & =_{(3)} &
  \sum_{\ba \in \gA^{\sV}}
        \mu\lb[\ba\rb] \cdot \barchi (\ba)  \\
 & = &
 \inn{\mu,\ \chi} 
\eeq

  (1) By definition of $\chi.\bg$ (equation (\ref{char.action})).

  (2) Because $\shift{\bg}:\gA^{\sU} \into \gA^{\bg.\sU}$ is an isomorphism. 

  (3) Because $\mu$ is locally $\dG-$invariant.

\proofof{``$\seilpmi$''} 
  If $\sV \subset \sU$ is finite, then for any $\ba \in \gA^\sV$,
then it is easy to verify that:

\[  \mu[\ba]  \ \ = \ \
 \marg{\sV}\mu [\ba] 
  \ \  =  \ \
 \sum_{\chi \in \dual{\sV}} \hmu_\chi \cdot \chi(\ba) \]

 The argument is then very similar to that of ``$\implies$''. 
\end{thmproof}

\subsection{Fourier Properties of Stationary Extensions}

  Suppose that $\sU \subset \sM$, and $\sV \subset \sU$ is a finite subset,
and suppose that $\chi \ \ := \ \ \Tensor_{v \in \sV} \ \chi_v$
  is some element of $\dual{\sU}$.  Then we can also think of $\chi$ as
an element of $\dual{\sM}$.  In other words, $\dual{\sU}$ embeds canonically
in $\dual{\sM}$.  We will ``abuse notation'', and identify elements of
$\dual{\sU}$ with their images in $\dual{\sM}$.  The following theorem
is a straightforward computation:

\Theorem{\label{harm.ext}}
{
  Let $\mu_\sU \in \Meas{\gA^\sU; \ \Cplx }$, and let $\mu \in
\Meas{\gA^\sM; \ \Cplx }$.  Then 
$\displaystyle  \statement{$\marg{\sU} \mu \ = \ \mu_\sU$} \ \ \iff  
  \ \ \statement{$\forall \chi \in \dual{\sU}, \ \ 
      \inn{\mu,\ \chi} \ = \ \inn{\mu_\sU, \ \chi} $}$. 
}

\qed

  Thus, we have reduced the Extension Problem to finding a measure
$\mu$ on $\Meas{\gA^\sM}$ whose Fourier coefficients agree with those
of $\mu_\sU$ on $\dual{\sU}$.  However, we can't just ``fill in'' the
remaining Fourier coefficients in an arbitrary way.  First of all,
we must produce something which is $\dG-$invariant.  Second of all,
we want to end up with a probability measure.

\Theorem{}
{
 Let $\mu_\sU \in \invMeas{\dG}{\gA^\sU}$, and let $\mu \in \Meas{\gA^\sM}$.
Then $\mu$ is a stationary extension of $\mu_\sU$ if and only if the
following two conditions are satisfied:

\bitem
  \item  For every $\chi \in \dual{\sU}$, and every $\bg \in \dG$, \ \ 
$ \inn{\mu, \ \chi.\bg} \ \ = \ \ \inn{\mu_\sU, \ \chi }$.

(This equation must be true even when $\chi.\bg$ is no longer in $\dual{\sU}$).

  \item  The Fourier coefficients of $\mu$ form a {\bf positive definite}
sequence.
\eitem
}
\begin{thmproof}
  The first condition follows from \thmpart{1} of Theorem \ref{stat.four}. 
Notice that, if more than one $\dG-$translate of $\chi$ lies inside
$\dual{\sU}$, then all of them will produce the same equation, by 
\thmpart{2} of Theorem \ref{stat.four} (since $\mu_\sU$ is locally 
$\dG-$invariant).

  The second condition is just the Bochner-Herglotz theorem to
guarantee that the measure $\mu$ is {\bf nonnegative}
\cite{Katznelson}.  This forces $\mu$ to be a probability measure,
because now $\mu[\gA^\sM] \ = \ \inn{\mu, \chr{} } \ = \ \inn{\mu_\sU,
      \chr{} } \ = \ \mu_\sU[\gA^\sU] \ = \ 1.$ (since $\mu_\sU$
      itself is a probability measure).
\end{thmproof}

\section{Extension on Finite Modules \label{finite} }

  Suppose that $\sM$ is a {\em finite} $\dG-$module, $\sU \subset \sM$, and
 and $\mu_\sU \in \invMeas{\dG}{\gA^\sU}$.  We will show that if
$\mu_\sU$ is ``sufficiently close'' to a product measure, then 
it is extendible.  More generally, we will show:

\Theorem{\label{thm.fin}}
{
 Let $\nu_\sU \in \invMeas{\dG}{\gA^\sU}$ be an extendible measure,
with an invariant extension $\nu$ such that $\supp{\nu} = \gA^\sM$.

  There exists an $\eps > 0$ so that, if $\mu_\sU \in \invMeas{\dG}{\gA^\sU}$
is any measure with $\norm{\mu_\sU - \nu_\sU}{v a r} \ < \ \eps$, 
then $\mu_\sU$ is also extendible.  This $\eps$ is of the form:

   \[  \eps \ \ = \ \ \frac{1}{ H(\sM) } \cdot \min_{\ba \in \gA^\sM} \nu[\ba]
                            \]

   ($\min_{\ba \in \gA^\sM} \nu[\ba] >0$ by hypothesis that  $\supp{\nu} = \gA^\sM$), where
$H(\sM)$ is a number determined by the $\dG-$module structure
of $\sM$, and which satisfies the following bounds:

\bitem
 \item (A) \ \ \ $H(\sM) \ \leq \ \card{\widehat{\gA^\sM}}$.

  \item (B) \ \ \ $H(\sM) \ \leq \  \card{\dG/\dH} \cdot \card{\widehat{\gA^\sU}}$.
\eitem

  where $\dH$ is the {\bf stabiliser} of $\sM$ in $\dG$:

  \[ \dH \ := \ \set{\bh \in \dG}{\forall m \in \sM, \ \bh \cdot m = m } \]

} 
\begin{thmproof}
  Define $\del_\sU := \mu_\sU - \nu_\sU$.  Thus, $\del_\sU$ is a real-valued measure.
Since $\mu_\sU$ and $\nu_\sU$ are locally $\dG-$invariant, $\del_\sU$ is also\footnote{
Cylinder subsets of $\gA^\sU$ can have negative $\del_\sU-$measures, but these
measures are still preserved under any shift which leaves the cylinder
set inside $\gA^\sU$.}.

  Next we will define $\del$, a real-valued, $\dG$-invariant measure
upon $\gA^\sM$, in terms of its Fourier coefficients.  For every $\chi \in
\dual{\sM}$, 

\bitem
  \item If there is some $\kappa \in \dual{\sU}$ and $\bg$ in $\dG$ so that
$\chi \ = \ \kappa.\bg$, then let 
$\hdel(\chi) \ \ := \ \ \widehat{\del_\sU}(\kappa)$.

  \item Otherwise, let $\hdel_\chi \ \ := \ \ 0$.
\eitem

  By \thmpart{2} of Theorem \ref{stat.four}, the definition of $\hdel_\chi$  
is independent of the choice of $\kappa$ and $\bg$, if more than one
choice is available.  By \thmpart{1} of the same theorem,  
the measure $\del$ is $\dG-$invariant.

\Claim{$\del$ is a real-valued measure.} 
\begin{claimproof}
  Since $\del_\sU$ is a real-valued measure, we know that,
for every $\chi \in \dual{\sU}, \ \ \widehat{\del_\sU} \lb( \barchi \rb) 
\ \ = \ \ 
\overline{\widehat{\del_\sU}( \chi)}$. 
It follows that, 
for every $\chi \in \dual{\sM}, \ \ \widehat{\del} \lb( \barchi \rb) 
\ \ = \ \ 
\overline{\widehat{\del}( \chi)}$, 
and from this, we conclude that $\del$ is also a real-valued measure.
\end{claimproof}

\claim{There is a number $H(\sM)$, determined by the $\dG-$module structure
of $\sM$, and satisfying inequalities (A) and (B), so that
$\norm{\del}{v a r} \ \ \leq \ \ H(\sM) \cdot \norm{\del_\sU}{v a r}$.
} 
\begin{claimproof}
  From elementary harmonic analysis \cite{Katznelson} , we know that:

\bitem
  \item $\norm{\widehat{\del_\sU}}{\oo} \ < \ \norm{\del_\sU}{v a r}$.

  \item $\norm{\del}{v a r} \ < \ \norm{\widehat{\del}}{1}$.
\eitem

 Hence, it suffices to show that $ \norm{\widehat{\del}}{1} \ < \
H(\sM) \cdot \norm{\widehat{\del_\sU}}{\oo}$, where $H(\sM)$ is the
aforementioned number.  To see inequality (A), notice that

\[ \norm{\widehat{\del}}{1} \ \leq \ \card{\sM} \cdot \norm{\widehat{\del}}{\oo}
                            \ = \ \card{\sM} \cdot \norm{\widehat{\del_\sU}}{\oo} \]

  where the second equality follows immediately from the definition of
$\hat{\del}$.

  Now for inequality (B).  For any $\chi \in
\widehat{\gA^\sU}$, let $\dG.\chi := \set{
\bg.\chi}{\bg \in \dG}$ be the {\bf orbit} of $\chi$ under the action of
$\dG$.  Then:

\beq
 \norm{\hat{\del}}{1} 
 & = &
 \sum_{\chi \in \widehat{\gA^\sM}} \lb| \hat{\del} (\chi) \rb| \\
 & = &
 \sum_{\chi \in \widehat{\gA^\sU}} \ \ \ \sum_{\xi \in \ \dG.\chi} \lb| \hat{\del} (\xi) \rb| \\
 & = &
 \sum_{\chi \in \widehat{\gA^\sU}} \ \ \ \sum_{\xi \in \ \dG.\chi} \lb| \widehat{\del_\sU} (\chi) \rb| \\ 
 & = &
 \sum_{\chi \in \widehat{\gA^\sU}} \card{\dG.\chi} \cdot \lb| \widehat{\del_\sU} (\chi) \rb| \\
\eeq

  But for any $\chi \in
\widehat{\gA^\sU}, \ \ \card{\dG.\chi} \ < \ \card{\dG / \dH}$.
So this expression is less than

\beq
 \sum_{\chi \in \widehat{\gA^\sU}} 
 \card{\dG/\dH} \cdot \lb| \widehat{\del_\sU} (\chi) \rb| 
&  = & 
  \card{\dG/\dH} \cdot \norm{ \widehat{\del_\sU}}{1} \\
& \leq &
  \card{\dG/\dH} \cdot \card{\widehat{\gA^\sU}} \cdot \norm{ \widehat{\del_\sU}}{\oo}
\eeq
\end{claimproof}
 
  Recall that $\nu$ is some invariant extension of $\nu_\sU$.  Define:

  \[ \mu \ \ :=  \ \ \nu + \del \]

\claim{$\mu$ is a nonnegative, $\dG-$invariant probability measure.}
\begin{claimproof}
  $\mu$ is a sum of two real-valued, $\dG-$invariant
 measures, and thus is also a real-valued, $\dG-$invariant measure.

Also, $\norm{\nu - \mu}{v a r} \ = \ \norm{\del}{v a r} \ <$ \
$H(\sM) \cdot \norm{\del_\sU}{v a r} \ = \ H(\sM) \cdot 
\norm{\nu_\sU - \mu_\sU}{v a r}$.  Thus, 

\beq
\lefteqn{ \lb(\maketall\norm{\nu_\sU - \mu_\sU}{v a r} <
\eps \ := \ \frac{1}{ H(\sM) } \cdot
\min_{\ba \in \gA^\sM} \ \ \nu[\ba] \rb)} \\
& \hspace{5em}\implies &
 \statement{ For every $\ba \ \in \gA^\sM, \ \ \mu[\ba] > 0$.} \\
\eeq

 It remains to show that $\mu[\gA^\sM] = 1$, or, equivalently,
that $\inn{\mu, \chr{}} \ = \ 1$.  Since  $\inn{\nu, \chr{}} \ = \ 1$,
this is equivalent to showing that  $\inn{\del, \chr{}} \ = \ 0$.
But $\inn{\del, \chr{}} \ = \ \inn{\del_\sU, \chr{}}$, and
$\inn{\del_\sU, \chr{}} \ = \ \inn{\nu_\sU, \chr{}} - \inn{\mu_\sU, \chr{}} \ = \ 0$.
\end{claimproof}

  Finally, we want to show that $\mu$ is an extension of $\mu_\sU$.
But
\[ \marg{\sU}{\mu} \ = \  \marg{\sU}{\nu} +  \marg{\sU}{\del}
                    \ = \ \nu_\sU + \del_\sU \ = \ \mu_\sU. \]
\end{thmproof}

  If $\rho$ is a probability measure on $\gA$, let $\rho^\sU$ be the
corresponding product measure on $\gA^\sU$.

\Corollary{}
{
 Let $\sM$ and $H(\sM)$ be as in the previous theorem.
Let $\rho$ be a probability measure on $\gA$ with full support, and
let
  \[ \eps \ \ := \ \ \frac{1}{H(\sM)}
          \lb( \min_{a \in \gA} \rho(a)\rb)^\card{\sM}  \]
  Let \ $\sU \subset \sM$.  If
$\mu \in \invMeas{\dG}{\gA^\sU}$, and $\norm{\mu - \rho^\sU}{v a r}
\ < \ \eps$, then $\mu$ is extendible.
}
\begin{thmproof}
 $\rho^\sU$ extends to the $\dG-$invariant probability measure 
$\rho^\sM$ on $\sM$, and 
$\displaystyle \min_{\ba \in \gA^\sM} \rho^\sM[\ba] \ \ = \ \ 
 \lb( \min_{a \in \gA} \rho(a)\rb)^\card{\sM}$. 
\end{thmproof}

\section{Envelopes:   Reduction to Smaller Modules
\label{envelope}}

  Suppose that $\sM$ and $\tlsM$ are $\dG-$modules, and that $\phi:\sM \into
\tlsM$ is a $\dG-$module {\bf homomorphism} ----that is, for
all $m \in \sM$ and $\bg \in \dG, \ \phi(\bg.m) \ = \ \bg.\phi(m)$.

  If $\tlba := \mtrx{\tla_\tlm}{\tlm \in \tlsM}{} \ \in \gA^\tlsM$, 
then define the element $\ba := \mtrx{a_m}{m \in \sM}{} \ \in \gA^\sM$,
by the formula:

  \beqn
\label{pullback.conf}
 \forall m \in \sM, \ \ a_m \ := \ \tla_{\phi(m)} 
  \eeqn

  This determines a function $\gA^\phi:\gA^\tlsM  \into  \gA^\sM$,
where $\gA^\phi(\tlba) \ := \  \ba $.

  If $\tlmu$ is a $\dG-$invariant measure on $\gA^\tlsM$, we define the
{\bf pullback} of $\tlmu$ through $\phi$ to be the measure:
 $\pb{\phi} \tlmu \ := \ (\gA^\phi)^* \mu$.  It is easily verified
that $\pb{\phi} \tlmu$ is a $\dG-$invariant measure on $\gA^\sM$.

  Given a $\dG-$module $\sM$, a subset $\sU
\subset \sM$ and a locally $\dG-$invariant
measure $\mu_\sU$ on $\gA^\sU$, we want to find a smaller $\dG-$module
$\tlsM$, a subset $\tlsU \subset \tlsM$, and a locally $\dG-$invariant
measure $\tlmu_\tlsU$ on $\gA^\tlsU$, such that, if we can extend $\tlmu_\tlsU$
to a $\dG-$invariant measure $\tlmu$ on $\gA^\tlsM$,  then 
$\mu := \pb{\phi} \tlmu$ is an extension of $\mu_\sU$

\Definition{Envelope}
{
  Let $\sM$ be a ~$\dG-$module, and $\sU \subset \sM$.

  An {\bf envelope} for $\sU$ is a~ $\dG-$module $\tlsM$, along with
a~ $\dG-$module homomorphism $\phi:\sM \into \tlsM$, such that

\bitem
   \item (E1) When restricted to $\sU$, the function $\phi$ is injective.

   \item (E2) If $\sV \subset \sU$, then for any $\tlbg \in \dG$ such that
   $\tlbg.\phi(\sV) \ \subset \ \phi(\sU)$, we can find some element $\bg \in \dG$
so that:

    \blist 
      \item $\bg.\sV \subset \sU$,
      \item For all $v \in \sV, \ \ \ \phi(\bg.v) \ = \ \tlbg . \phi(v).$
(Thus, $\phi(\bg.\sV) \ =  \ \tlbg.\phi(\sV)$.)
    \elist
 
\eitem
}

\subparagraph{Example:}{\em Envelopes in a Lattice}

  Suppose $\dG \ = \ \sM \ = \ \Zahl^D$, and let $\sU \subset \Zahl^D$ be
finite, and small enough that it fits into a box of dimensions
$N_1 \x N_2 \x \ldots \x N_D$.  We will indicate the action of $\Latt$ on
itself with the ``+'' symbol. 

Consider the $\Zahl^D-$module: 

   \[ \tlsM \ \ := \ \ \frac{\Zahl}{2 N_1 \Zahl} \ \x \ 
      \frac{\Zahl}{2 N_D \Zahl} \ \x \ \ldots \x
      \frac{\Zahl}{2 N_D \Zahl} \]

  and let $ \phi: \sM \into \tlsM$ be the  $\Zahl^D-$module homomorphism: 

  \[ \phi(n_1,\ldots,n_D) \ := \
       \lb( n_1 + \frac{\Zahl}{2 N_1 \Zahl}, \ \
           n_2 + \frac{\Zahl}{2 N_2 \Zahl}, \  \ \ldots , \ \
           n_D + \frac{\Zahl}{2 N_D \Zahl} \rb) \]

  Then $(\tlsM, \phi)$ is an envelope for $\sU$.

\subparagraph{Remark:}  In this example, the module

  \[ \tlsM \ \ := \ \ \frac{\Zahl}{ N_1 \Zahl} \ \x \ 
      \frac{\Zahl}{ N_2 \Zahl} \ \x \ \ldots \x
      \frac{\Zahl}{ N_D \Zahl} \]

  with the quotient map $\phi:\sM \into \tlsM$ would {\em
not} necessarily have worked as an envelope for $\sU$.  To see this,
suppose that

  \[  \sU \ \ := \ \ \CC{1..N_1} \x \{1\} \x \{1\} \x \ldots \x \{1\} \]

 and let $ \sV \ \ := \ \ \{ \bv_1, \bv_2 \},$ where $\bv_1 \ := \
(1,1,\ldots,1)$, while $\bv_2 \ := \ (2,1,1,\ldots,1)$.  Let $\tlbg \ :=
\ (N_1\!-\!1,\ 0,\ 0,\ \ldots,\ 0) \ \in \ \dG$.  Then note that

\[ \tlbg + \phi(\bv_1)
 \ \  = \ \ 
  \phi( \tlbg + \bv_1 ) 
 \ \  = \ \ 
 \phi(N_1,1,1,\ldots,1)
 \ \  = \ \ 
 \phi(\bv_3), \]

  where $\bv_3 := (N_1,1,1,\ldots,1)$, while

\[ \tlbg + \phi(\bv_2) 
 \ \  = \ \ 
  \phi( \tlbg + \bv_2 ) 
 \ \  = \ \  
\phi(N_1+1,1,1,\ldots,1)
 \ \  = \ \  
\phi(1,1,1,\ldots,1) 
 \ \  = \ \ 
\phi(\bv_1) \]

  Now, there is no element $\bg \in \dG$ so that 
$\bg + \sV \ \ = \ \ \lb\{ \bv_1, \bv_3 \rb\}. $
Thus, although $\tlbg + \phi(\sV) \subset \phi(\sU)$, we {\em cannot}
find some $\bg \in \dG$ so that $\bg + \sV \subset \sU$ and $\phi(\bg + \sV) \ = \ \tlbg + \phi(\sV)$.

\pause

\Proposition{\label{thm.env}}
{
  Let $\sM$ be a~ $\dG-$module, and $\sU \subset \sM$.  Let $\phi:M \into \tlsM$
be an envelope for $\sU$, and $\tlsU := \phi(\sU)$.

\blist
  \item  For any probability measure $\mu_\sU$ on $\gA^\sU$, there is a unique
probability measure $\tlmu_\tlsU$ on $\gA^\tlsU$ so that
$ \mu_\sU \ \ = \ \ \pb{\phi} \tlmu_\tlsU.$

  \item  If $\mu_\sU$ is locally $\dG-$invariant, then so is $\tlmu_\tlsU$.

  \item If $\tlmu$ is an extension of $\tlmu_\tlsU$ to a~ $\dG-$invariant
probability measure on $\gA^\tlsM$, then $\nu \ := \ \pb{\phi} \tlmu$ is an
extension of $\mu_\sU$ to a~ $\dG-$invariant probability measure on
$\gA^\sM$,
\elist
}

\begin{thmproof}

\proofof{\thmpart{1}}  By hypothesis,  $\phi_| : \sU \into \tlsU$ is injective.  Let
 $\psi: \tlsU \into \sU$ be the inverse map, and define $\tlmu_\tlsU
:= \pb{\psi} \mu_\sU$.  Thus,
$\mu_\sU \ = \ \pb{\phi} \tlmu_\tlsU$.  Since $\phi_{|\sU}$ is injective,
the measure $\tlmu_\tlsU$ is the unique one satisfying this equation.

\proofof{\thmpart{2}}
 Let $\tlsV \subset \tlsU$, and $\tlbc \in \gA^\tlsV$.
 Suppose $\tlbg \in \dG$ is
such that $\tlbg.\tlsV \subset \tlsU$ as well.  We want to show:

   \[ \tlmu_\tlsU \lb[ \shift{\tlbg}\,\tlbc \rb] \ \ = \ \ \tlmu_\tlsU \lb[ \tlbc \rb] \]

  Let $\sV \ := \ \psi (\tlsV) \ \subset \ \sU$, and let $\bc \ := \ \gA^\phi(\tlbc)$, where $ \gA^\phi:\gA^\tlsV \into \gA^\sV$ is as defined by
equation (\ref{pullback.conf}) near the
beginning of \S\ref{envelope}.  Thus,  if $\tlbc \ = \ \seq{\tlc_v}{v
\in \tlsV}$, then
 $\bc \ = \ \seq{c_v}{v \in \sV}$, where, for all 
$v \in \sV, \ \ c_v := \ \tlc_{\phi(v)}$.

    Let $\bC$ be the cylinder set in $\gA^\sU$ associated to $\bc$
(and likewise, $\tlbC$ for $\tlbc$).  Thus, $\tlbC \ = \
\gA^\psi(\bC).$ Since $\tlsM$ is an envelope, there is a $\bg \in \dG$
satisfying condition (E2).  By (E2)(1), $\shift{\bg}\,\bC$ is also a cylinder
set in $\gA^\sU$, and since $\mu_\sU$ is locally $\dG-$invariant,
$\mu_\sU \lb[ \shift{\bg}\,\bC \rb] \ = \ \mu_\sU[\bC]$.

\Claim { $\gA^\psi(\shift{\bg}\,\bC) \ = \ \shift{\tlbg}\,\tlbC$}

\begin{claimproof}
Let $\tlba := \mtrx{\tla_\tlu}{\tlu \in \tlsU}{} \ \in \gA^\tlsU$,
and suppose that $\tlba = \gA^\psi(\ba)$, where  
$\ba := \mtrx{a_u}{u \in \sU}{} \ \in \gA^\sU$.  Then
$\statement{$\tlba \in \gA^\psi(\shift{\bg}\,\bC)$}
\iff$ 
 $\statement{$\ba \in \shift{\bg}\,\bC$} 
\iff$ 
$\statement{$\forall v \in \sV, \ a_{\lb(\bg.v\rb)} \ = \ c_v$}
\iff_{\!\!\!(1)}$ \\
$ \statement{$\forall v \in \sV, \ \ 
  \tla_{\lb(\tlbg. \phi(v)\rb)} \ = \ \tlc_{\phi(v)}$ }
\iff$  $\statement{$\tlba \in \shift{\tlbg}\,\tlbC $}.$

(1) Because, for all $v \in \sV$, \ 
$ \tla_{\lb(\tlbg. \phi(v)\rb)} \ = \ \tla_{\phi(\bg.v)}
\ = \ a_{\bg.v}$, \ \ and $ c_v = \tlc_{\phi(v)}$.
\end{claimproof}

  Thus, $\tlmu_\tlsU \lb[ \shift{\tlbg}\,\tlbC \rb]
 \ = \ 
 \tlmu_\tlsU \lb[ \gA^\psi(\shift{\bg}\,\bC) \rb] 
 \ \  = \ \ 
 \mu_\sU \lb[ \shift{\bg}\,\bC \rb] 
 \ = \ 
 \mu_\sU \lb[ \bC \rb]
 \ =  \ 
 \tlmu_\tlsU \lb[ \tlbC \rb]. $

\proofof{\thmpart{3}}  This is straightforward.
\end{thmproof}

\section{Embedding of Locally Stationary Measures \label{embed}}

  Suppose that $(X,\sX,\nu)$ is a probability space, and $T$ is
a $\nu$-preserving action of $\Latt$ upon $X$.  Let 
$\sP:X \into \gA$ be a measurable function (ie. a
$\gA$-labelled, {\bf measurable partition} of $X$), and
let $\sP^\Latt:X \into \gA^\Latt$ be the map
$x \mapsto \mtrx{\sP\lb(T^\bn(x) \rb)}{\bn \in \Latt}{}$.
The projection of $\mu$ through $\sP^\Latt$
is then a stationary probability measure on $\gA^\Latt$,
called the {\bf stochastic process induced by $\sP$ and $T$}.
Call this measure $\eta$.

  Suppose that $\sU \subset \Latt$, and $\mu_\sU \in
\statMeas{\gA^\sU}$.  The map $\sP$ is an {\bf embedding} of $\mu_\sU$
in the system $(X,\sX,\nu; \ T)$ if $\marg{\sU}{\eta} = \mu_\sU$.
When can $\mu_\sU$ be thus embedded?

\Theorem{\label{thm.embed}}
{
  Suppose that $\sU \subset \Latt$ is finite, and that $\mu_\sU$ lies
in the {\bf interior} of $\extMeas{\gA^\sU}$.  
Suppose that $(X,\sX,\nu; \ T)$ is {\bf ergodic}.  Then $\mu_\sU$
can be embedded in $(X,\sX,\nu; \ T)$.
}
\begin{thmproof}
  We will first show how to construct an ``approximate'' embedding for
$\mu_\sU$.  The approximation method involves a certain degree of error,
which can be exactly characterized and then compensated for. 

  Suppose $U \in \Natur$, so that $\sU \subset \sB(U)$.
Let $\mu \in \statMeas{\gA^\Latt}$ be an extension of $\mu_\sU$.  
Then for any $N > 0, \ \ \mu_{\sB(N)} := \marg{\sB(N)} \mu$ is a
locally stationary probability measure on $\gA^{\sB(N)}$.  Also, 
if $\sU_0 \subset \sB(N)$ is any translation of $\sU$, then 
$\marg{\sU_0}{\mu_{\sB(N)}} \ = \ \mu_{\sU_0}$,  where $\mu_{\sU_0}$ is
the obvious ``translation'' of $\mu_\sU$ to the domain $\sU_0$.

  The {\bf Rokhlin Tower Lemma} for $\Latt$-actions says that,
for any $\eps > 0$ and $N \in \Natur$, there is a subset $R \in \sX$
so that the disjoint union:

 \[ \Disj_{\bn \in \sB(N+U)} T^\bn (R) \]

  has measure greater than $1-\eps$.

   Let $x \in X$ be a generic point for $R$, and suppose we look at the
``name'' of $x$ with respect to the partition $\{R, \ X \setminus R\}$:
for all $\bn \in \Latt$, colour the point $\bn$ ``black'' if $T^\bn x \in
R$, and ``white'' otherwise.  Let $\sR \subset \Latt$ be the set of
``black'' points.  The Rokhlin Tower condition is equivalent to saying
that the union:

  \[ \Disj_{\br \in \sR} \lb( \maketall \sB(N+U)+\br \rb) \]

 is disjoint, and has \Cesaro density greater than $1 - \eps$ in $\Latt$.

  To define a measurable function $\sP:X \into \gA$, we will provide a
scheme to determine its value at every point in the $\Latt-$orbit
of $x$, in terms of the $\{R, \ X \setminus R\}-$name of $x$ (this is 
sometimes called ``colouring the name of $x$'').  The scheme well-defines
the values of $\sP$ on the orbit of every generic point in $X$ ---thus,
it defines $\sP$ almost everywhere on $X$.

  Defining the value of $\sP$ on the $\Latt-$orbit of $x$ is equivalent to
defining a function $\bp: \Latt \into \gA$ ---in other words, a
configuration.  Do this as follows:  Let $\phi:\sR \into \gA^{\sB(N)}$
be some function so that, for each $\ba \in \gA^{\sB(N)}$, the
\Cesaro density of the subset $\phi^{-1}(\ba)$ inside $\sR$ is equal to
$\mu_{\sB(N)} [\ba]$ (since the set $\sR$ itself has a well-defined
\Cesaro density, such a function can always be constructed).  
For each $\bu \in \sR$, let $\bp_{\sB(N)+\bu} = \phi(\bu)$.  This
immediately defines $\bp$ on ``most'' of $\Latt$.  Now, fix some $\ga \in
\gA$, and label {\em all} remaining points in $\Latt$ with the symbol $\ga$.

   The function $\sP$ induces a stationary probability measure $\eta$
on $\gA^\Latt$. \ \ $\eta_\sU := \marg{\sU} \eta$ is ``close'' to
$\mu_\sU$, but slightly ``enriched'' in words that contain big blocks
of the ``$\ga$'' symbol, while impoverished in words that don't.  If
we fix $\eps > 0$ and $N \in \Natur$, then
$ \eta_\sU \ = \ F_{\eps,N} [\mu_\sU]$,
where $F_{\eps,N}:\extMeas{\gA^\sU} \into \extMeas{\gA^\sU}$ is an
{\bf affine function}.

 So, if we want to actually produce the measure $\mu_\sU$ as an {\em
outcome} of this procedure, we must find some $\nu_\sU \in
\extMeas{\gA^\sU}$, so that $\mu_\sU = F_{\eps,N} [\nu_\sU]$.  In other
words, in order to use this construction to build an embedding of
$\mu_\sU$ within $X$, we must find some $N$ and $\eps$ so that $\mu_\sU
\in I_{N,\eps} := F_{N,\eps}\lb( \extMeas{\gA^\sU} \rb)$.

\claim{For any $\del > 0$, there exist $\eps$ and $N$ so that
$\Lbsg[I_{\eps,N}] \ \geq \
 (1-\del)\cdot\Lbsg \lb[\extMeas{\gA^\sU}\rb]$,
 where $\Lbsg$ is the Lebesgue measure.}
\begin{claimproof} 
   $F_{\eps,N}$ is affine, and thus, differentiable with a constant
derivative, $D_{\eps,N}$.  For any $\delta_1 > 0$, we can find a small
enough $\eps$ and large enough $N$ that, for every $\mu_\sU \in
\extMeas{\gA^\sU}$, \ $\norm{ F_{\eps,N} [\mu_\sU] \ - \ \mu_\sU}{v a
r} < \del_1$.  Thus, for any $\del_2 > 0$, we can make $\del_1$ small enough so
that
$\norm{D_{\eps,N} - \Id{}}{\oo} < \del_2$ (where $\norm{\cdot}{\oo}$ is
the operator norm).  Thus, for any $\del$, we can in turn make $\del_2$
small enough that the determinant of $D_{\eps,N}$ is within $\del$ of $1$.
Thus, for large enough $N$ and small enough $\eps, \ \ F_{\eps,N}:
\extMeas{\gA^\sU} \into I_{\eps,N}$ is a diffeomorphism, and, if $\Lbsg$
is the Lebesgue measure, then $\Lbsg[I_{\eps,N}] \ \geq \ (1-\del)\cdot
\Lbsg \lb[\extMeas{\gA^\sU}\rb]$.  
\end{claimproof}

\claim{For any $\mu$ in the interior of $\extMeas{\gA^\sU}$ there exist $\eps$
and $N$ so that $\mu \in I_{\eps,N}$.}
\begin{claimproof}
Identify $\Meas{\gA^\sU; \ \Real}$ with $\Real^{\gA^\sU}$, endowed with an
inner product.  $I_{\eps,N}$ is convex, so if $\mu \in \extMeas{\gA^\sU}
\setminus I_{\eps,N}$, then there is some unit vector $\bv \in
\Real^{\gA^\sU}$, so that $I_{\eps,N} \subset
\set{ \bw \in \Real^{\gA^\sU}}{ \inn{ \bw - \mu, \ \bv} < 0 }$.
Fix $\mu$, and regard $m_\bv$ as a function of $\bv$.
The set \\  $\set{ \bw \in  \Real^{\gA^\sU}}{ \inn{ \bw - \mu, \ \bv} \geq 0 }
\intsct \extMeas{\gA^\sU}$ has nontrivial Lebesgue measure
$m_\bv \cdot \Lbsg \lb[\extMeas{\gA^\sU}\rb]$, for some $m_\bv > 0$.
Since the unit sphere in $\Real^{\gA^\sU}$ is compact, there is some
$M > 0$ so that $m_\bv \geq M$ for all $\bv$ in the sphere.

  Let $\del < M$, and, by Claim $1$, find $\eps$ and $N$ so that
$\Lbsg[I_{\eps,N}] \geq (1-\del)\cdot
\Lbsg \lb[\extMeas{\gA^\sU}\rb]$.  Then we have $M \cdot \Lbsg \lb[\extMeas{\gA^\sU}\rb]$ \ $> \ \del \cdot \Lbsg \lb[\extMeas{\gA^\sU}\rb]$ \ $> \
\Lbsg \lb[\extMeas{\gA^\sU} \setminus I_{\eps,N} \rb]$ \\ $\geq \
\Lbsg\lb[ \set{ \bw \in  \Real^{\gA^\sU}}{ \inn{ \bw - \mu, \ \bv} > 0 } 
\intsct \extMeas{\gA^\sU} \rb]$ \ $> \ M \cdot \Lbsg \lb[\extMeas{\gA^\sU}\rb]$,
a contradiction.
\end{claimproof}

  We conclude that any point $\mu$ in the interior of $\extMeas{\gA^\sU}$ 
is in $I_{\eps,N}$ for some $\eps$ and $N$, and thus, can be
``embedded'' in the system $(X,\sX,\mu; \ T)$ via the aforementioned
construction.
\end{thmproof}

\section{(quasi)Periodic, Ergodic, and Mixing Extensions
\label{sect.per}}

\subsection{Periodic Probability Measures
\label{per.prob.meas}}

  If $\dP \subset \Natur^D$, then a configuration
$\ba \in \gA^\Latt$ is called {\bf $\dP$-periodic} if, for all $n \in \Latt$
and $p \in \dP$, 
\ \ \ $a_{n+p} =  a_{n}$.  If $\inn{\dP}$ is the sublattice generated
by $\dP$, and $\tlsM := \Latt/\inn{\dP}$, with 
$\Latt$ acting upon $\tlsM$ by translation, then $\tlsM$ is
$\Latt-$module.  The quotient map $\phi:\Latt \into \tlsM$ is a
homomorphism of $\Latt-$modules.  Configuration $\ba$ is
$\dP$-periodic if and only if $\ba \ = \ \gA^\phi \tlba$, for some
word $\tlba \in \gA^\tlsM$ (in the notation of Section
\ref{envelope}).

  In general, if $\sM$ is a $\dG-$module, $\tlsM$ is another $\dG-$module,
and $\phi:\sM \into \tlsM$ is a $\dG-$module homomorphism, then
we will say that an element $\ba \in \gA^\sM$ is {\bf $\tlsM$-periodic}
if  $\ba \ = \ \gA^\phi [\tlba]$, for some $\tlba \in \gA^\tlsM$.

  If $\mu$ is a $\dG-$invariant measure on $\gA^\sM$, then 
$\mu$ is {\bf $\tlsM$-periodic} if the elements of the space
$\lb(\gA^\sM, \mu \rb)$ are $\mu-$almost surely $\tlsM-$periodic.
This is the case if and only if there is a $\dG-$invariant
measure $\tlmu$ on $\gA^\tlsM$, such that 
$\mu \ = \ \pb{\phi} [\tlmu]$.

\subsection{Periodic Extensions}

  Suppose that $\sU \subset \sM$, and $\mu_\sU$ is a locally $\dG-$invariant
measure upon $\gA^\sU$.  Can we extend $\mu_\sU$ to a periodic measure
on $\gA^\sM$?

\Theorem{\label{ext.per}}
{
  Suppose that $\tlsM$ is a {\em finite} $\dG-$module, a quotient of
$\sM$ via the map $\phi:\sM \into \tlsM$, and an {\bf envelope} for $\sU$.
Let $H(\tlsM)$ be the constant described in Theorem \ref{thm.fin}

  Let $\tlnu$ be a $\dG-$invariant measure on $\gA^\tlsM$, with full
support, and let 

  \[ \eps \ := \
    \frac{1}{H(\sM)} \cdot \min_{\tlba \in \gA^\tlsU} \tlnu[\tlba] \]

  Let $\nu \ = \ \lb(\gA^\phi\rb)^* \tlnu$, and let $\nu_\sU := \marg{\sU} \nu$.
If $\mu_\sU$ is any locally $\dG-$invariant measure on $\gA^\sU$ so that
$ \norm{\mu_\sU - \nu_\sU}{v a r} \ < \ \eps$,
then $\mu_\sU$ can be extended to a $\dG-$invariant, $\tlsM-$periodic
probability measure on $\gA^\sM$.
}
\begin{thmproof}
Let $\tlsU := \phi(\sU) \subset \tlsM$.  By \thmpart{1} of Theorem
\ref{thm.env}, the measure $\tlmu_\tlsU := \pb{(\phi^{-1})} \mu_\sU$
is a locally $\dG-$invariant measure on $\gA^\tlsU$.  Further, if
$\tlnu_\tlsU := \marg{\tlsU} \tlnu$, then $\norm{\tlmu_\tlsU -
\tlnu_\tlsU}{v a r} < \eps$.  Since $\tlsM$ is finite, we can apply
Theorem \ref{thm.fin}, and extend $\tlmu_\tlsU$ to a $\dG-$invariant measure,
$\tlmu$, on all of $\gA^\tlsM$.

  Now, define $\mu := \ \pb{\phi} [\tlmu]$.  Then $\mu$ is 
a $\tlsM-$periodic, $\dG-$invariant measure by construction, and also,
$\marg{\sU} \mu \ = \ \mu_\sU$.
\end{thmproof}

\Corollary{}
{
 The set $\extMeas{\gA^\sU}$ has nontrivial {\bf interior} in the
space $\Meas{\gA^\sU; \ \Real}$, and the set of $\tlsM$-periodically
extendible measures has nontrivial interior within $\extMeas{\gA^\sU}$.
}
\begin{thmproof}
  Let $\rho$ be any probability measure on $\gA$ with full support,
and let $\mu_\sU := \rho^\sU$ be the product measure on $\gA^\sU$.  In
the notation of Theorem \ref{ext.per}, $\rho^\tlsM$ is a
$\dG$-invariant extension of $\rho^\tlsU$, with full support, and
induces a $\tlsM$-periodic extension of $\mu_\sU$ to $\gA^\sM$.  By
Theorem \ref{ext.per}, all measures in an open ball around $\mu_\sU$
also have $\tlsM$-periodic extensions.
\end{thmproof}

\Corollary{}
{
  Suppose $\sU \subset \Latt$ is finite, and fits inside a box of size
$Q_1 \x Q_2 \x \ldots \x Q_D$.  Suppose that $\bP := (P_1,\ldots,P_D)$, where
$P_1 \geq 2 Q_1, \   P_2 \geq 2 Q_2, \ \ldots ,  P_2 \geq 2 Q_2$, and
let $\nu$ be a $\bP-$periodic, stationary probability measure on $\gA^\Latt$.
Let $\nu_\sU := \marg{\sU} \nu$.

  There is an $\eps > 0$ (a function of $\bP$ and $\nu$), so that,
if $\mu_\sU$ is any locally stationary probability measure on $\gA^\sU$ 
within $\eps$ of $\nu_\sU$ in total variation norm, then $\mu_\sU$
has a $\bP-$periodic extension.
}
\qed

  For any $\bP := (P_1,\ldots,P_D)$, let 
$\invMeas{\bP}{\gA^\Latt}$ denote the set of $\bP-$periodic, stationary
processes.

 If $\sU \subset \Latt$, then let $\invMeas{\bP}{\gA^\sU}$ denote the
set of {\bf $\bP$-periodically-extendible measures}:  those elements of
$\statMeas{\gA^\sU}$ having an extension that is $\bP-$periodic.
The following facts are not difficult to verify:

\bitem
  \item $\invMeas{\bP}{\gA^\sU}$ is a closed, convex set.
 
  \item If $\mu \in \invMeas{\bP}{\gA^\sU}$ and 
$\nu \in \invMeas{\bQ}{\gA^\sU}$, then any convex combination of $\mu$
and $\nu$ is inside $\invMeas{\bR}{\gA^\sU}$, where, for each $d \in \CC{1...D}, \ \ 
R_d$ is the {\em lowest common multiple} of $P_d$ and $Q_d$. 
\eitem

   Let $\invMeas{per}{\gA^\sU}$ be the set of all locally stationary
measures possessing a periodic extension of any periodicity.
It follows that $\invMeas{per}{\gA^\sU}$ is also a convex set.

\subsection{Essentially Aperiodic measures
\label{essentially.aperiodic}}

  Not every extendible measure has a periodic extension.  This
follows from the existence of {\bf essentially aperiodic tile systems}
---that is, sets of tiles which can tile the plane, but {\em only}
in an aperiodic fashion.  In \cite{RMR}, Raphael Robinson exhibits
a collection of six ``notched'' square tiles, which, along with their
4 rotations, will tile the plane, but {\em only} in an aperiodic fashion.
We can code these six tiles as six $3 \x 3$ matrices in the alphabet
$\gA := \{ 0, \fa,\fA,\fb,\fB,\fc,\fC \}$

\[ \begin{array}{ccc}
\frame{$\begin{array}{ccc}
\fA & \fC & \fA \\
\fB & 0 & \fd \\
\fA & \fB & \fA \\
\end{array}$}
&
\frame{$\begin{array}{ccc}
\fa & \fc & \fa\\
\fc & 0 &  \fc \\
\fa & \fC & \fa \\
\end{array}$}
&
\frame{$\begin{array}{ccc}
\fa & \fb & \fa\\
\fc & 0 &  \fc \\
\fa & \fB & \fa \\
\end{array}$}
\\
\frame{$\begin{array}{ccc}
\fa & \fC & \fa\\
\fB & 0 &  \fC \\
\fa & \fB & \fa \\
\end{array}$}
&
\frame{$\begin{array}{ccc}
\fa & \fb & \fa\\
\fc & 0 &  \fc \\
\fa & \fb & \fa \\
\end{array}$}
&
\frame{$\begin{array}{ccc}
\fa & \fb & \fa\\
\fb & 0 &  \fb \\
\fa & \fB & \fa \\
\end{array}$}
\\
\end{array}\]

Each tile has a ``0'' symbol in its center, surrounded by four ``corners''
and four ``edges''.  The tiles must be put together so that these corners
and edges ``match'' according to the following {\bf mapping rules:}

\bitem
  \item  ``$\fb$'' edges must be matched to  ``$\fB$'' edges.
  \item  ``$\fc$'' edges must be matched to  ``$\fC$'' edges.
  \item  Where four tiles meet, exactly three corners must be
       of type ``$\fa$'', and one of type ``$\fA$''.
\eitem

  These matching rules can be encoded as a subshift of finite type on
the alphabet $\gA$, defined by some subset $\gR \subset \gA^\sU$, where
$\sU := \CC{1..3}^2$.   Any configuration in $\inn{\gR}$
corresponds to some Robinson tiling.  Now let $\mu$ be a stationary
probability measure on $\inn{\gR}$, and let $\mu_\sU := \marg{\sU} \mu$.
Then $\mu_\sU$ is a locally stationary measure, and $\supp{\mu_\sU} \ = \ \gR$.

  We claim that $\mu_\sU$ is ``essentially aperiodic''.  To see this,
suppose that $\nu$ was any extension of $\mu_\sU$.  Then $\supp{\nu} \
\subset \inn{\gR}$, and thus, almost every configuration in the
probability space $(\gA^{\Zahl^2},\nu)$ is aperiodic.

\subsection{Essentially Periodic Measures}

  At the opposite extreme are {\bf essentially periodic} measures:
locally stationary measures which {\em only} have periodic extensions.

  For example, let $\gA := \{0,1\}$ and $\sU \ := \ \CC{1...9} \x \{0,1\}$, and
let $\gB \subset \gA^\sU$ be the set:

\[ \lb\{ 
  \begin{array}{c}
  \frame{${\scriptsize \begin{array}{ccccccccc}
            0&0&0&0&0&0&0&0&0\\
            1&1&1&1&1&1&1&1&0\\  
           \end{array}}$}, \\
  \vdots\\ 
  \frame{${\scriptsize \begin{array}{ccccccccc}
            0&0&0&0&0&1&0&0&0\\  
            0&0&0&0&0&0&1&1&0\\
           \end{array}}$},\\
  \frame{${\scriptsize \begin{array}{ccccccccc}
            0&0&0&0&0&0&1&1&0\\  
            0&0&0&0&0&0&1&0&0\\
           \end{array}}$},\\ 
  \frame{${\scriptsize \begin{array}{ccccccccc}
            0&0&0&0&0&0&1&0&0\\  
            0&0&0&0&0&0&0&1&0\\
           \end{array}}$},\\ 

    \frame{${\scriptsize \begin{array}{ccccccccc}
            0&0&0&0&0&0&0&1&0\\  
            0&0&0&0&0&0&0&0&0\\
           \end{array}}$} \\ 
  \end{array}
   \rb\} \]

   In other words, all blocks are of the form 
\frame{${\scriptsize \begin{array}{cc}w_1 & 0 \\ w_0 & 0\end{array}}$} 
where $w_0$ and $w_1$ are successive $8-$bit binary numbers.
Let $\gW \subset \gA^\sU$ be the set containing all elements of $\gB$ and
all their horizontal cyclic permutations.  $\gB$ defines a subshift
of finite type, which contains only the orbit of
a single, periodic configuration, having horizontal periodicity
$9$, and vertical periodicity $256$.  Call this configuration $\ba$

  If $\mu_\sU$ is the measure on $\gA^\sU$ assigning equal mass
to each of the $2304$ elements of $\gB$, then $\mu_\sU$ has only one
stationary extension:  the measure $\mu$ which assigns equal mass
to each of the $2304$ distinct translates of $\ba$.  Thus,
$\mu_\sU$ is {\bf essentially periodic}, with period $256 \x 9$.

  Note that the periodicity $256 \x 9$ is much larger than $2 \x 9$, which
was the size of the initial domain $\sU$.  Indeed, as this argument
makes clear, the periodicity of essentially periodic measure can be made
to grow exponentially with the size of the initial domain.

\subsection{Ergodic Extensions \label{sect.erg.extend}}

  A stationary probability measure $\mu$ on $\gA^\Latt$ is called 
{\bf ergodic} if any measurable subset $\bU \subset \gA^\Latt$ which
is invariant under all shifts has $\mu-$measure either zero or one.
The set of ergodic measures on $\gA^\Latt$, which we denote by 
``$\ergMeas{\gA^\Latt}$'', is exactly the set of {\em extremal
points} of $\statMeas{\gA^\Latt}$ (see \cite{Furstenberg} or \cite{Walters})
Hence, every stationary measure can be approximated arbitrarily
well as a convex combination of ergodic measures.

  If $\sU \subset \Latt$, and $\mu_\sU \in \statMeas{\gA^\sU}$, then 
we say $\mu$ is {\bf ergodically extendible} if it can be extended
to an ergodic measure on $\gA^\Latt$.  The set of ergodically
extendible measures will be written as ``$\ergMeas{\gA^\sU}$''. 
Since the map $\pr{\sU}^*:
\Meas{\gA^\Latt} \into \Meas{\gA^\sU}$ is linear, any extremal
point of $\extMeas{\gA^\sU}$ has a $\pr{\sU}^*-$preimage which is
extremal in $\statMeas{\gA^\Latt}$.  As a consequence, every
extremal point of $\extMeas{\gA^\sU}$ is in $\ergMeas{\gA^\sU}$.
Hence, every extendible measure on $\gA^\sU$ can be approximated arbitrarily
well as a convex combination of ergodically extendible measures.

\breath

 We will see in Section \ref{mix.extend} that, in fact, ``almost all''
extendible measures are ergodically extendible.  However, not every
extendible measure is.  To see this, suppose that $\sU \subset \Latt$
is some finite domain, let $\gA$ and $\gB$ be two {\em disjoint}
alphabets, and suppose that $\mu_\sU \in \statMeas{\gA^\sU}$ and
$\nu_\sU \in \statMeas{\gB^\sU}$ are two extendible probability
measures.  Let $\eta_\sU := \frac{1}{2} \mu_\sU + \frac{1}{2}
\nu_\sU$.  Then $\eta_\sU$ is also extendible, and any extension of
$\eta_\sU$ is of the form $\eta := \frac{1}{2} \mu + \frac{1}{2} \nu$,
where $\mu$ and $\nu$ extend $\mu_\sU$ and $\nu_\sU$, respectively.
$\eta$ can never be ergodic: $\ \gA^\Latt$ and $\gB^\Latt$ are
disjoint, shift-invariant subsets of $\lb(\gA \disj \gB\rb)^\Latt$,
each having $\eta-$measure $\frac{1}{2}$.

\Proposition{}
{
Let $\sU \subset \Latt$ be finite.
\blist
  \item  Every ergodically extendible measure on $\gA^\sU$ is a limit point of
periodically extendible measures.
  \item  $\invMeas{per}{\gA^\sU}$ is a dense, convex subset of 
$\extMeas{\gA^\sU}$.  
  \item 

 $\invMeas{per}{\gA^\sU}$ contains the entire {\bf interior} of
$\extMeas{\gA^\sU}$.
\elist
 }
\begin{thmproof}  \thmpart{2} follows immediately from \thmpart{1}, and
the fact that $\invMeas{per}{\gA^\sU}$ is convex, and the fact that
$\extMeas{\gA^\sU}$ is the convex closure of $\ergMeas{\gA^\sU}$.

\proofof{\thmpart{3}} This follows from \thmpart{2}, and the fact that, if $C$
a dense, convex subset of a $D$-dimensional convex set $K$, then $C$
contains $\interior{K}$.  To see this, let $x \in \interior{K}$, and
let $B$ be an open ball around $x$ inside of $\interior{K}$.  Let $S$
be the boundary of $B$, and let $s_1,\ldots,s_D$ be $D$ equidistant
points in $S$, so that their convex closure, $co\{s_1,\ldots,s_D\}$ is
a regular $D$-simplex containing the centre-point $x$.
 
  Since $C$ is dense in $K$, \ $C \intsct B$ is dense in $B$.  
Thus, find elements $c_1,\ldots,c_D \in C$ so that, for
all $d \in \CC{1..D}$, \ $c_d$ is ``very close'' to $s_d$.  Then
$co\{c_1,\ldots,c_D\} \subset C$ is a $D$-simplex ``very close'' to
$co\{s_1,\ldots,s_D\}$,  and therefor contains $x$.

\proofof{\thmpart{1}}  Let $\mu_\sU \in
\ergMeas{\gA^\sU}$, and let $\mu$ be an ergodic extension of
$\mu_\sU$.  Let $\ba \in \gA^\Latt$ be a {\em generic configuration}
for $\mu$: in other words, for any finite subset $\sV \subset \Latt$
and configuration $\bb \in \gA^\sV$,

  \[ \mu[\bb] \ \ = \ \ \lim_{N \goto \oo} \Freq{\bb}{\ba}{\sB(N)} \]

  where $\sB(N) := \CO{0...N}^D$ is the $D-$dimensional cube of
side length $N$, and 

\beq
 \Freq{\bb}{\ba}{\sB(N)} & := & 
  \frac{\btxt{ \# of times ``$\bb$'' appears inside $\ba_{\sB(N)}$ \etxt}} 
     {N^D}  \\
& = &
  \frac{1}{N^D} \sum_{\bn \in \sB(N)} \chr{}\{\ba_{\sV+\bn} = \bb\}
\eeq
   
  Such generic configurations exist, by the Birkhoff Ergodic Theorem.

  In particular, for any $\eps > 0$, we can find a large enough $N$ so that,
for {\em all} $\bb \in \gA^\sU$, 

  \[  \lb| \maketall
   \mu[\bb] - \Freq{\bb}{\ba}{\sB(N)} \rb| \ < \ \frac{\eps}{2} \]

  Suppose that all of $\sU$ fits inside a cube of side length $U$.
Assume that $N$ is so large that the $U-$thick boundary of $\sB(N)$ is
``relatively small'':

  \[ \frac{\card{\sB(N)} - \card{\sB(N-U)}}{\card{\sB(N)}} \ < \
  \frac{\eps}{2} \] 

  Now, identify $\sB(N)$ with $\sN := (\Zahl/N) \dirsum \ldots \dirsum (\Zahl/N)$,
and treat $\ba_{\sB(N)}$ as an element of $\gA^\sN$.  Then this configuration,
along with its $N^D$ periodic translations on $\gA^\sN$, defines a 
stationary measure on $\gA^\sN$, which, in turn, defines an $\sN-$periodic,
stationary measure on $\gA^\Latt$.  Call this measure $\nu$, and then
let $\nu_\sU := \marg{\sU} \nu$. It is straightforward to verify that 

  \[ \norm{ \nu_\sU \ - \ \mu_\sU  } \ < \ \eps \]

  and of course, by construction, $\nu_\sU \in \invMeas{per}{\gA^\sU}$.
\end{thmproof}

\subsection{Mixing, Weak Mixing, and Quasiperiodicity
\label{mix.extend}}

  A stationary probability measure $\mu$ on $\gA^\Latt$ is called 
{\bf weakly mixing} if the stochastic process $(\gA^\Latt \x \gA^\Latt,\ 
\mu \tensor \mu)$ is ergodic.  $\mu$ is called {\bf mixing} if,
for any measurable $A,B \subset \gA^\Latt$ of nonzero measure,
any any sequence $\seq{\bn_k}{k \in \Natur} \subset \Latt$ tending
to infinity, $\lim_{k \goto \oo} \mu\lb[ A \intsct \shift{\bn_k} B\rb]
= \mu[A] \cdot \mu[B].$  A function $\phi \in \bL^2(\gA^\Latt, \mu)$ 
is an {\bf eigenfunction} of the system  $(\gA^\Latt, \mu)$ 
if there is a group homomorphism $\chi:\Latt \into \Torus{1}$
such that, for all $\bn \in \Latt$, \ \ $\shift{\bn} (\phi)
= \chi(\bn) \cdot \phi$.  The system is called {\bf quasiperiodic}
if $\bL^2(\gA^\Latt, \mu)$ has an orthonormal basis of eigenfunctions.

  All of these concepts can be defined for any measure-preserving
$\Latt-$action on a probability space $(X,\sX,\mu)$.  Mixing
implies weak mixing implies ergodicity, but weak mixing and
quasiperiodicity are mutually exclusive.
Furthermore, all of these properties are {\bf inheritable}
through morphisms.  If $(X,\sX,\mu; \ T)$ and $(\hX,\hsX,\hmu; \ \hT)$ are
two measure-preserving $\Latt-$actions, then a {\bf morphism} between
the systems is a measure-preserving surjection $\Psi:X \into \hX$ so
that, for all $\bn \in \Latt, \  \Psi \circ T^{\bn} = \hT^{\bn} \circ
\Psi$.  If $\Psi$ is such a morphism, and $(X,\sX,\mu; \ T)$ is
ergodic (respectively: weakly mixing, mixing, or quasiperiodic),
then so is $(\hX,\hsX,\hmu; \ \hT)$.

  In particular, let $F: X \into \gA$ be a measurable function, so
that $F$ and $T$ together {\bf induce} a stationary stochastic process
on $\gA^\Latt$, having measure $\hmu$ (see Section \ref{embed}).  If
$\hX := \supp{\hmu} \subset \gA^\Latt$ and $\hT:=
\shift{}$, then the map $F^\Latt:X \into \gA^\Latt$ is a morphism.
Thus, if $(X,\sX,\mu; \ T)$ possesses any of the aforementioned
inheritable properties, so does the process $(\gA^\Latt, \hmu)$.

\Theorem{}
{ 
 Suppose that $\sU \subset \Latt$ is finite, and that $\mu_\sU$ is 
   in the {\bf interior} of $\extMeas{\gA^\sU}$.  Then
   $\mu_\sU$ can be extended to a stationary process $\mu$ which
is any of:  ergodic, mixing, weakly mixing, or quasiperiodic.
}
\begin{thmproof}
  The argument is the same in all four cases.  First, find a 
system $(X,\sX,\nu; \ T)$ which is ergodic, and which also has the
property in question (for the first three, this is trivial;  for
the fourth, it is sufficient to know that ergodic, quasiperiodic
systems exist).  Next, use Theorem \ref{thm.embed} to 
{\bf embed} $\mu_\sU$ within the desired process.  Let $\mu \in
\statMeas{\gA^\Latt}$ be the stochastic process generated by
this embedding.  Then $\mu$ itself has the desired property.
\end{thmproof}

The same argument works for any other ``inheritable'' property of
dynamical systems.  The interpretion: knowledge of the local marginal
$\mu_\sU$ tells you basically nothing about the asymptotic dynamical
properties of the process $\mu$.

\section{Decidability Questions \label{decide}}

  In section \ref{finite.type}, we showed:

\bquote
{\em It is formally undecidable whether, for a given subset $\gW
\subset \gA^\sU$, the set $\extMeas{\gW}$ is nonempty.}
\equote

  This raises the question of whether the Extension Problem itself is
formally decidable.

  Let $\RReal$ be the set of all {\bf recursively computable} ({\bf r.c}) real
numbers: that is, real numbers whose decimal expansion can be
generated by some Turing Machine \cite{HopcroftUllman}.  $\RReal$ is a
countable field, containing all rational and real-algebraic numbers.
Let $\RMeas{\gA^\sM; \ \Real}$ be the set of {\bf r.c.}, real-valued
measures: those such that, if $\sV \subset \sM$ is finite, and $\ba
\in
\gA^\sV$, then the measure of $\ba$ is an element of $\RReal$.
(Of course, some ``exotic'' measurable subsets of $\gA^\Latt$ may
have non {\bf r.c.} measures). 
$\RMeas{\gA^\sM; \ \Real}$ is a vector space over the field $\RReal$.

  Let $\RinvMeas{\dG}{\gA^\sM}$ be the set of $\dG-$invariant probability
measures, etc.  Clearly, when we ask about the ``formal decidability'' of
the Extension Problem, what we are {\em really} referring to is the
Extension Problem for {\bf r.c.} measures:

\bquote
{\em If $\sU \subset \sM$, and $\mu_\sU \in \RinvMeas{\dG}{\gA^\sU}$, 
is $\mu$ extendible to a $\dG-$invariant measure on $\gA^\sU$?}
\equote

  Note that we do not require the extension itself to be {\bf r.c}.
If a recursive decision procedure explicitly {\em constructs} an extension,
then this extension will be {\bf r.c.} by nature.  However, it is
conceivable that some recursive decision procedure might exist which
demonstrates the existence of an extension by ``nonconstructive''
means.  It is conceivable that, although we can recursively
decide that $\mu_\sU$ is extendible, no {\bf r.c} extension exists.

  A subset $\bS \subset \RinvMeas{\dG}{\gA^\sU}$ is called
{\bf recursively decidable} ({\bf r.d}) if there is a Turing machine $\dM$,
so that, when given any $\mu \in \RinvMeas{\dG}{\gA^\sU}$ as input,
$\dM$ halts after some {\em finite} number of steps, and outputs
either ``yes'' or ``no'', depending upon whether or not $\mu$ is
an element of $\bS$.

   A subset $\bS \subset \RinvMeas{\dG}{\gA^\sU}$ is called {\bf
recursively enumerable} ({\bf r.e}) if there is a Turing machine
$\dM$, so that, when given any integer $n \in \Natur$ as input, $\dM$
halts after a finite number of steps, and produces as output some 
measure $F_\dM[n] \in
\bS$, and so that the function $F_\dM:\Natur \into \bS$ instantiated by
$\dM$ is {\em surjective.}  In other words, $\dM$ provides a mechanism to
systematically ``list'' all elements of $\bS$.

  Equivalently, $\bS \subset \RinvMeas{\dG}{\gA^\sU}$ is 
{\bf recursively enumerable} if there is a Turing machine $\dM$,
so that, when given any $\mu \in \RinvMeas{\dG}{\gA^\sU}$ as input,
$\dM$ halts after some finite number of steps {\em unless}
$\mu$ is {\em not} in $\bS$, in which case $\dM$ never halts.

   The following facts are easy to verify: Any {\bf r.d}
set is {\bf r.e.}, but the converse is not true.  However,
if both $\bS$ and its complement are {\bf r.e.}, then $\bS$ is
{\bf r.d}.  Finally, although a countable union of
{\bf r.d} sets is not necessarily itself {\bf r.d}, it is still
{\bf r.e.} \cite{HopcroftUllman}.

\Theorem{}
{
 Let $\sU \subset \Latt$ be a finite subset.  Then

\blist
  \item  For any $\bP \in \Natur^D, \ \ \RinvMeas{\bP}{\gA^\sU}$ is
{\bf r.d}.

  \item   $\RinvMeas{per}{\gA^\sU}$ is {\bf r.e.}.

  \item  $\RstatMeas{\gA^\sU} \setminus  \RextMeas{\gA^\sU}$
is {\bf r.e.}.
\elist
}
\begin{thmproof}

\proofof{\thmpart{1}} 
  If $\mu_\sU \in \RstatMeas{\gA^\sU}$, we want to know whether
the set
$S := \set{ \mu \in \invMeas{\bP}{\gA^\Latt}}{\marg{\sU} \mu = \mu_\sU}$ 
is nonempty.

  Suppose $\bP := (P_1,\ldots,P_D)$.  Let $\tlsM := (\Zahl/P_1) \dirsum \ldots
\dirsum (\Zahl/P_D)$, and suppose that $\sU$ maps bijectively into the
subset $\tlsU \subset \sM$ via the quotient map from $\Latt \into \tlsM$.
Let $\tlmu_\tlsU \in \statMeas{\gA^\tlsU}$ be the projected image of $\mu_\sU$.

  The vector space of $\bP-$periodic, signed measures on $\gA^\Latt$ is
linearly isomorphic to the finite dimensional vector space
$\Meas{\gA^\tlsM; \ \Real}$.  The image of $S$ under this isomorphism
is the affine set 
\[\tlS := \set{\mu \in \Meas{\gA^\tlsM; \ \Real}}{\mu \
\btxt{a stationary probability measure, and \etxt} \marg{\tlsU} \mu =
\tlmu_\tlsU}.\] 
  $\tlS$  is the solution set of a finite system of linear
equations and linear inequalities in $\mu$:

\bitem
  \item $\mu\lb[\gA^\tlsM\rb] = 1$.
  \item  For all $n \in \Latt, \ \ \shift{n}_* \mu = \mu$. 
  \item  $\marg{\tlsU} \mu = \tlmu_\tlsU$.
  \item  For all $\ba \in \gA^\sM, \ \ \mu[\ba] \geq 0$. 
\eitem

  Thus, it is {\bf r.d} whether $\tlS$ is nonempty, and thus,
whether $\mu_\sU$ has a $\bP-$periodic extension.

\proofof{\thmpart{2}} 
  $\RinvMeas{per}{\gA^\sU}$ is a countable union of recursively
decidable sets, and thus, {\bf r.e.}.

\proofof{\thmpart{3}} 
  Suppose $\mu \in \Meas{\gA^\Latt; \ \Cplx}$ has Fourier transform
$\mtrx{\hmu_\chi}{\chi \in \dual{\Latt}}{}$, let $\sV \subset \Latt$
be a finite subset, and let $\ba \in \gA^\sV$.  It is easy
to verify:

  \[ \mu[\ba] \ = \ \sum_{\chi \in \dual{\sV}} 
 \hmu_\chi \cdot \overline{\chi(\ba)} \]

  Thus, if $\mu_\sU \in \statMeas{\gA^\sU}$, then by Theorems 
\ref{stat.four} and \ref{harm.ext}, $\mu$ is an extension of $\mu_\sU$ if and only if:

\bitem
   \item For all $\chi \in \dual{\sU}, \ \ \hmu_\chi \ = \
   \inn{\mu_\sU, \chi}$

   \item For all  $\bn \in \Latt$, and all $\chi \in \dual{\Latt}$, if
$\xi := \chi \circ \shift{\bn}$ then $\hmu_\chi \ = \ \hmu_\xi$.

   \item For all finite $\sV \subset \Latt$ and $\ba \in \gA^\sV, \
\ \ \sum_{\chi \in \dual{\sV}}  \hmu_\chi \cdot \overline{\chi(\ba)} \ >\  0$.
\eitem

  Thus, an extension for $\mu_\sU$ is equivalent to a set of Fourier
coefficients satisfying a countable collection of linear equations and
inequalities.

  For all $N \in \Natur$, let $\sB(N) := \CC{0...N}^D$, and let $\Xi_N :=
\dual{\sB(N)}$.  If $N$ is large enough that
$\sU \subset \sB(N)$, then we can start by trying to define all the
Fourier Coefficients in the set $\set{\mu_\chi}{\chi \in \Xi_{N}}$.  The
three sets of linear constraints listed above now become a finite
system of linear equations and inequalities ---if the solution set is
nonempty, call it $S_N$. 

\Claim{Suppose that, for all $N \in \Natur$, the set $S_N$ is nonempty.
Then $\mu_\sU$ is extendible.}
\begin{claimproof}
$S_N$ is a compact subset of the
finite dimensional vector space $\Cplx^{\Xi_N}$.  Furthermore, if
$S_{N+1}$ is also nonempty, then any vector in $S_{N+1}$, when projected to 
$\Cplx^{\Xi_N}$, determines an element in $S_N$.  Call this projection map
$\pr{N}$.

  Fix $N$, and, for all $M > N$, let $\tlS^M_N := \pr{N} \circ \pr{N+1}
\circ \ldots \circ \pr{M-1} (S_M)$, a nonempty compact subset of
$S_N$.  Also, $\tlS^{M+1}_N \supset \tlS^{M+2}_N
\supset \tlS^{M+3}_N \supset \ldots$. \  Thus, $\tlS_N := \intsct_{M > N}
\tlS^M_N$ is a nonempty compact subset.  Further,
$\pr{N} \lb(\tlS_{N+1}\rb) = \tlS_N$.  Thus, any element of $\tlS_N$ can be
``extended'' to an element of $\tlS_{N+1}$, which can then be
``extended'' to $\tlS_{N+1}$, etc.  

  Pick any element $\hmu_N \in \tlS_N$, and inductively extend it in
this fashion, producing $\hmu_M \in \tlS_M$, for every $M > N$.  
Once this is done, the collection of vectors $\seq{ \hmu_M }{M > N}$
defines a single element $\hmu \in \Cplx^{\dual{\Latt}}$.  
\ $\hmu$ is the Fourier transform of some measure $\mu$, and by construction,
$\mu$ is a stationary probability measure, and an extension of $\mu_\sU$.
\end{claimproof}

  Hence, if $\mu_\sU$ is {\em not} extendible, then, by contradiction,
there must be some $N \in \Natur$ so that $S_N$ is empty.
Since $S_N$ is the solution set of a finite system of linear equations
and inequalities, it is {\bf r.d} whether $S_N$ is
empty.  

  Hence, by successively checking the nonemptiness of $S_N$ for
each $N \in \Natur$, we have a recursive procedure which will
halt {\em if} $\mu_\sU$ is {\em not} extendible, and tell us so.
(If $\mu_\sU$ {\em is} extendible, however, the procedure will never
halt).  Thus we can recursively enumerate the elements of
$\RstatMeas{\gA^\sU} \setminus  \RextMeas{\gA^\sU}$.
\end{thmproof}

\Theorem{}
{
Let $\sU \subset \Latt$ be finite.  The set
$\RextMeas{\gA^\sU}$ is {\em not} {\bf r.e.}.
}
\begin{thmproof}
Recall that, if $\gT \subset \gA^\sU$, then  $\inn{\gT}$ is the
associated {\bf subshift of finite type} (see Section \ref{finite.type}).
Let $\bN := \set{\gT}{\inn{\gT} \ \btxt{is not trivial\etxt}}$,
and let $\bT := \set{\gT}{\inn{\gT} \ \btxt{is trivial\etxt}}$.
Recall that $\bN$ is not {\bf r.d} (see \cite{RMR}, \cite{BerDomino},
or \cite{KitS}).

\claim{Suppose $\gT \subset \gA^\sU$.   If $\gT \in \bT$, then there
is some $N \in \Natur$ so that no configuration in $\gA^{\sB(N)}$ is
$\gT-$admissable}
\begin{claimproof}
  Suppose that, for every $N \in \Natur$, there was a configuration
$\ba^{[N]} \in \gA^{\sB(N)}$ that was $\gT-$admissable ---that is:
for all $n \in \sB(N)$, if $n+\sU \subset \sB(N)$, then $\ba^{[N]}_{n+\sU} \in
\gT$.  Extend $\ba^{[N]}$ to an element of $\gA^\Latt$ by filling all the
remaining entries in some arbitrary fashion ---call the extended
configuration $\bb^{[N]}$ 

  Since $\gA^\Latt$ is compact, the sequence $\mtrx{\bb^{[N]}}{N \in \Natur}{}$
has a convergent subsequence ---call it $\mtrx{\bb^{[N_k]}}{k \in \Natur}{}$
----which converges to some limit $\bb \in \gA^\Latt$.

  For any $M \in \Natur$, there is some $K \in \Natur$ so that, for all $k > K,
\ \ \bb^{[N_k]}_{\sB(M)} = \bb_{\sB(M)}$.  Hence, the central ``$\sB(M)-$block''
of $\bb$ is $\gT-$admissable.  This is true for every $M$; 
we conclude that $\bb$ is $\gT-$admissable.  Thus, the set $\inn{\gT}$ is
nonempty, since it contains $\bb$.
\end{claimproof}

\claim{The set $\bT$ is {\bf r.e.}.}
\begin{claimproof}
  Fix $\gT \subset \gA^\sU$.
For any finite $N$, it is {\bf r.d} whether
or not $\gA^{\sB(N)}$ contains a $\gT-$admissable configuration (there are
only a finite number of cases to check).
Suppose we perform this procedure for every $N \in \Natur$.  By 
Claim $1$, if $\gT \in \bT$, then we will eventually find an $N$
where no $\gT-$admissable configuration exists.  Thus, we have
a procedure which will halt if $\gT \in \bT$, and tell us so.
\end{claimproof}

  As a consequence, since $\bN$ is {\em not} {\bf r.d},
we conclude that $\bN$ is not even {\bf r.e.}.

\claim{Suppose that $\RextMeas{\gA^\sU}$ was {\bf r.e.}.
Then $\bN$ is also {\bf r.e.}.}
\begin{claimproof}
Clearly, $\bN = \set{ \gT \subset \gA^\sU}{ \btxt{for some \etxt} \ \mu \in 
\RextMeas{\gA^\sU}, \ \ \supp{\mu} = \gT}$.  Hence, any recursive
procedure for enumerating the elements of $\RextMeas{\gA^\sU}$
would also provide a means for enumerating the elements of $\bN$.
\end{claimproof}

  By contradiction, $\RextMeas{\gA^\sU}$ cannot be {\bf r.e.}.
\end{thmproof}

\section{Conclusion
\label{conclusion}}

  Although $\extMeas{\gA^\sU}$ itself is not recursively denumerable,
both its complement and topological interior are (Section
\ref{decide}).  As yet, however, no efficient procedure exists for
determining when a locally stationary measure is extendible.  So far
the only substantive result in this direction is Theorem
\ref{ext.per}, which says, loosely, that if $\mu_\sU$ is
``sufficiently close'' to a periodically extendible measure with full
support, then $\mu_\sU$ itself is periodically extendible.

  The existence of mixing, ergodic, etc. extensions is
well-characterized in Section \ref{mix.extend}.  However, as yet, no
useful work has been done characterizing the {\bf entropy} of these
extensions.  In particular, we might ask: given that $\mu_\sU$ is
extendible, what do the {\bf maximal-entropy} extensions of $\mu_\sU$
look like?  Is the maximal-entropy extension {\em unique}?  Does it
possess some kind of ``Markov'' property, analogous to the Markov
Extension in $\Zahl$?  Perhaps it is some kind of {\bf Markov Random
Field} \cite{Rozanov}.  Indeed, in general, what would a ``Markov
extension'' of a locally stationary measure look like, if anything?
In the nonprobabilistic, purely symbolic setting, the construction
analogous to a Markov extension is a $\Zahl^D$-subshift of finite
type, but these are still poorly understood.  Even topological Markov
shifts ---the simplest subshifts of finite type ---do not generalize
easily to higher dimensions
\cite{MarkleyPaul1}.  The maximal entropy measures for such subshifts
have been studied in \cite{MarkleyPaul2}; \ perhaps similar techniques
can be applied to maximal-entropy extensions of probability measures.

\paragraph*{Acknowledgements:}

  I would like to thank Andres del Junco, Jeremy Quastel, and Reem
Yassawi for their advice and suggestions.

\bibliographystyle{plain}
\bibliography{../../bibliography}

\hrulefill

\begin{verbatim}
  Marcus Pivato
  Department of Mathematics,
  University of Houston

  Email: pivato@math.toronto.edu
\end{verbatim}
\end{document}